\documentclass[a4wide]{article}
\usepackage[utf8]{inputenc}
\usepackage[a4paper, total={6.5in, 9.in}]{geometry}

\usepackage{hyperref}
\usepackage[auth-lg]{authblk}
\usepackage{amsfonts}
\usepackage{babel}[english]
\usepackage{bm}
\usepackage{mathtools}
\usepackage{eufrak}
\usepackage{enumerate}
\usepackage{algorithm,algorithmic}
\usepackage{caption}
\usepackage{subcaption}
\usepackage{xcolor}
\usepackage{diagbox}
\usepackage{bigints}
\usepackage{relsize}
\usepackage{comment}

\newtheorem{remark}{Remark}
\usepackage{array}
\newcolumntype{x}[1]{>{\centering\arraybackslash\hspace{0pt}}p{#1}}

\begin{document}
	
\newcommand{\WidestEntry}{$\left[13,8,555 \right]$}
\newcommand{\SetToWidest}[1]{\makebox[\widthof{\WidestEntry}]{#1}}

\title{Deep-HyROMnet: A deep learning-based operator approximation for hyper-reduction of nonlinear parametrized PDEs}
\author[1]{Ludovica Cicci}
\author[1]{Stefania Fresca}
\author[1]{Andrea Manzoni}
\setlength{\affilsep}{1em}
\renewcommand\Authsep{, }
\affil[1]{MOX-Dipartimento di Matematica, Politecnico di Milano, P.zza Leonardo da Vinci 32, 20133 Milano, Italy\footnote{\texttt{{ludovica.cicci,stefania.fresca,andrea1.manzoni}@polimi.it}}}

\date{}
\renewcommand\Affilfont{\small}
 
\maketitle

\begin{abstract}

To speed-up the solution to parametrized differential problems, reduced order models (ROMs) have been developed over the years, including projection-based ROMs such as the reduced-basis (RB) method, deep learning-based ROMs, as well as surrogate models obtained via a machine learning approach.  
Thanks to its physics-based structure, ensured by the use of a Galerkin projection of the full order model (FOM) onto a linear low-dimensional subspace,  RB methods yield approximations that fulfill the physical problem at hand. However, to make the assembling of a ROM independent of the FOM dimension, intrusive and expensive hyper-reduction stages are usually required, such as the discrete empirical interpolation method (DEIM), thus making this strategy less feasible for problems characterized by (high-order polynomial or nonpolynomial) nonlinearities.  
To  overcome this bot\-tle\-neck, we propose a novel strategy for learning nonlinear ROM operators using deep neural networks (DNNs). The resulting hyper-reduced order model enhanced by deep neural networks, to which we refer to as Deep-HyROMnet, is then a physics-based model, still relying on the RB method approach, however employing a DNN architecture to approximate reduced residual vectors and Jacobian matrices once a Galerkin projection has been performed. Numerical results dealing with fast simulations in nonlinear structural mechanics show that Deep-HyROMnets are orders of magnitude faster than POD-Galerkin-DEIM ROMs, keeping the same level of accuracy.

\end{abstract}


\section{Introduction}

Mathematical models involving 
partial differential equations (PDEs) depending on a set of parameters are ubiquitous in applied sciences and engineering. These input parameters are defined to characterize, e.g., material properties, loads, boundary/initial conditions, source terms, or geometrical features. High-fidelity simulations based on full-order models (FOMs), like the finite element method (FEM), entail huge computational costs in terms of CPU time and memory, if a large number of degrees of freedom is required. Complexity is amplified whenever interested in going beyond a single direct simulation, such as in the multi-query contexts of optimization, parameter estimation and uncertainty quantification.  
To face these problems, several strategies to build reduced order models (ROMs) have been developed over the years, aiming at computing reliable solutions to parametrized PDEs at a greatly reduced cost.  

A large class of ROMs relies on a projection-based approach, which aims at approximating the unknown state quantities as a linear superimposition of basis functions; these latter then span a subspace which the governing equations are projected onto 
\cite{benner2015survey, benner2017model}. Among these, the reduced basis (RB) method \cite{quarteroni2016reduced, hesthaven2016certified} is a powerful and widely used technique, characterized by a splitting of the reduction procedure into an expensive, parameter-independent offline phase (however performed once and for all) and an efficient, parameter-dependent online phase.  Its efficiency mainly relies on two crucial assumptions:
\begin{enumerate}
\item the solution manifold is low-dimensional, so that the FOM solutions can be approximated as a linear combination of few reduced modes with a small error; 
\item  the online stage is completely independent of the high-fidelity dimension \cite{farhat20205}. 
\end{enumerate}

Assumption 1 concerns the approximability of the solution set and is associated with the slow decay of the  Kolmogorov $N$-width \cite{pinkus2012n}.  However, for physical phenomena characterized by a slow $N$-width decay, such as those featuring  coherent structures that propagate over time \cite{fresca2020deep}, the manifold spanned by all the possible solutions is not of small dimension, so that ROMs relying on   linear (global) subspaces might be inefficient. Alternative strategies to overcome this bottleneck can be, e.g., local RB methods \cite{amsallem2012nonlinear, pagani2018numerical, vlachas2021local}, or nonlinear approximation techniques, mainly based on deep learning architectures, see, e.g., \cite{lee2020model, kim2020efficient, fresca2021comprehensive, fresca2021pod, franco2021deep}.

Assumption 2 is automatically verified in linear, affinely parametrized problems \cite{quarteroni2016reduced}, but cannot be fulfilled when dealing with nonlinear problems, as the online assembling of the reduced operators requires to reconstruct the high-fidelity ones. To overcome this issue, a further level of approximation, or {\em hyper-reduction}, must be introduced. State-of-the-art methods, such as the empirical interpolation method (EIM) \cite{barrault2004empirical}, the discrete empirical interpolation method (DEIM) \cite{chaturantabut2010nonlinear}, its variant matrix DEIM \cite{negri2015efficient}, the missing point estimation \cite{astrid2008missing} and the Gauss-Newton with approximated tensors (GNAT) \cite{carlberg2011efficient}, aim at recovering an affine expansion of the nonlinear operators by computing only a few entries of the nonlinear terms. EIM, DEIM and GNAT can be seen as {\em approximate-then-project} techniques, since operator approximation is performed at the level of FOM quantities, prior to the projection stage. On the other hand, {\em project-then-approximate} strategies have also been introduced, aiming at approximating directly ROM operators, such as the reduced nonlinear residual and its Jacobian. An option in this sense is represented by the so-called Energy Conserving Sampling and Weighting (ECSW) technique \cite{farhat2015structure}.  See. e.g., \cite{farhat20205} for a detailed review. \\

Although extensively applied in many applications, spanning from fluid flows models to cardiac mechanics \cite{drohmann2012reduced, amsallem2012nonlinear, tiso2013discrete, radermacher2016pod, ghavamian2017pod, bonomi2017reduced}, these strategies are code-intrusive and, more importantly, might impact on the overall efficiency of the ROM approximation in complex applications. Very often, when dealing with highly nonlinear problems expensive hyper-reduction strategies are required if aiming at preserving the physical constraints at the ROM level, that is, if  ROMs are built consistently with the FOM through a projection-based strategy. For instance, a large number of DEIM basis vectors are required to ensure the convergence of the reduced Newton systems arising from the linearization of the nonlinear hyper-ROM when dealing with highly nonlinear elastodynamics problems \cite{cicci2021cardiacDEIM}, even if few basis functions are required to approximate the state solution in a low-dimensional subspace.  An alternative formulation of DEIM in a finite element (FE) framework, known as unassembled DEIM \cite{tiso2013modified}, has been proposed to preserve the sparsity of the problem, while in \cite{peherstorfer2014localized} a localized DEIM selecting smaller local subspace for the approximation of the nonlinear term is presented. 

Semi-intrusive strategies, avoiding the construction of a ROM through a Galerkin projection, have been recently proposed exploiting surrogate models to determine the RB approximation. For instance, neural networks (NNs) or Gaussian process (GP) regression can be employed to learn the map between the input parameters and the reduced-basis expansion coefficients in a non-intrusive way \cite{hesthaven2018non,guo2018reduced,guo2019data,swischuk2019projection}. An approximation of the nonlinear terms arising in projection-based ROMs is obtained in \cite{gao2020non} through deep NNs (DNNs) that exploit the projection of FOM solutions.  
NNs have also  also been recently applied in the context of {\em operator inference} for (parametrized) differential equations, combining ideas from classical model reduction with data-driven learning. For instance, the design of NNs able to accurately represent linear/nonlinear operators, mapping input functions to output functions, has been proposed recently in \cite{lu2021learning}; based on the universal ap\-pro\-xi\-ma\-tion theorem of operators \cite{chen1995universal}, a general deep learning framework, called DeepONet, has been introduced to learn continuous operators, such as solution operators of PDEs, using DNNs; see also \cite{wang2021learning}. In \cite{peherstorfer2016data} a non-intrusive projection-based ROM for parametrized time-dependent PDEs  including low-order polynomial nonlinear terms is considered, inferring an approximation of the reduced operators directly from data of the FOM. Finally, the obtained low-dimensional system is solved -- in this case, the learning task consist in the solution to a least squares problem; see also \cite{benner2020operator}. Projection-based ROMs and machine learning have been fused in \cite{qian2019transform} aiming at the approximation of linear and quadratic ROM operators, focusing on the solution to a large class of fluid dynamics applications. Similarly, in \cite{bai2021non} a non-intrusive technique, exploiting machine learning regression algorithms, is proposed for the approximation of ROM operators related to projection-based methods for the solution of parametrized PDEs. Finally, \cite{bhattacharya2021model} combines principal component analysis-based model reduction with a NNs for approximation, in a purely data-driven fashion, of infinite-dimensional solution maps, such as the solution operator for time-dependent PDEs. \\

In this paper, we develop a novel semi-intrusive, deep learning-enhanced hyper-reduced order modeling strategy, which hereon we refer to as Deep-HyROMnet, by leveraging a Galerkin-RB method for solution dimensionality reduction and DNNs to perform hyper-reduction. Since the efficiency of the nonlinear ROM hinges upon the cost-effective approximation of the projections of the (discrete) reduced residual operator and its Jacobian (when an implicit numerical scheme is employed), the key idea is to overcome the computational bottleneck associated with classical, intrusive hyper-reduction techniques, e.g. DEIM, by relying on DNNs to approximate otherwise expensive reduced nonlinear operators at a greatly reduced cost. Unlike data-driven-based methods, for which the predicted output is not guaranteed to satisfy the underlying PDE, our method is physics-based, as it computes the ROM solution by actually solving the reduced nonlinear systems by means of Newton method, thus exploit the physics of the problem. A further benefit of the method proposed lies on the fact that the inputs given to the NN are low-dimensional arrays, so that the overwhelming training times and costs that may be required by even moderately large FOM dimensions can be avoided. We point out that Deep-HyROMnet aims at efficiently approximate the nonlinear operators given by the composition of the reduced solution operator, that maps the input parameter vector and time to the corresponding ROM solution, and the reduced residual/Jacobian operator, that maps the ROM solution to the reduced residual/Jacobian evaluated on the ROM solution. To the best of our knowledge, this is the first method of its kind. We apply the novel methodology to the solution of problems in nonlinear solid mechanics, with particular focus on parametrized nonlinear elastodynamics and complex (e.g., exponential nonlinear) constitutive relations of the material undergoing large de\-for\-ma\-tions, showing that Deep-HyROMnet outperforms the DEIM-based ROM in terms of computational speed-up in the online stage, still achieving accurate results.  \\

The paper is structured as follows. We recall the formulation of the RB method for nonlinear unsteady parametrized PDEs in Section~\ref{sec:RBM}, relying on POD for the construction of the reduced subspace and on DEIM as hyper-reduction technique. Deep-HyROMnet and the DNN architecture employed to perform reduced operator approximation are detailed in Section~\ref{sec:Deep-HyROMnet}. The numerical performances of the method are assessed in Section~\ref{sec:tests} on several benchmark problems related with nonlinear elastodynamics. Finally, conclusions and future perspective are reported in Section~\ref{sec:conclusion}.


\section{Projection-based ROMs: the reduced basis method}\label{sec:RBM}

Our goal is to pursue an efficient solution to nonlinear unsteady PDE problems  depending on a set of input parameters, which can be written in abstract form as follows: given an input parameter vector $\bm\mu\in\mathcal{P}$, $\forall t\in(0,T]$, find   $\mathbf{u}(t;\bm \mu)\in V$ such that
\begin{equation} \label{eq:motion}
R(\mathbf{u}(t;\bm \mu),t;\bm\mu)=0 \quad \text{in } V',
\end{equation}
where the parameter space $\mathcal{P}\subset\mathbb{R}^P$ is a compact set and $H^1(\Omega)^m \subseteq V \subseteq H_0^1(\Omega)^m$ is a suitable functional space, depending on the boundary conditions at hand, whereas $\Omega \subset \mathbb{R}^d$ is a bounded domain in $d$ dimensions, $d = 1, 2, 3$. In particular, we are interested in vector problems ($m=3$) set in $d=3$ dimensions. The role of the parameter vector $\bm\mu$ depends on the particular application at hand; in the case of nonlinear elastodynamics, $\bm\mu$ is related to the coefficients of the constitutive relation, the material properties and the boundary conditions. 

By performing discretization in space and time, we end up with a fully-discrete nonlinear system 
\begin{equation}\label{eq:residual}
\mathbf{R}(\mathbf{u}_h^n(\bm\mu),t^n;\bm\mu) =  \mathbf{0} \quad \text{in } \mathbb{R}^{N_h},
\end{equation}
at each time step $t^n$,  $n=1,\dots,N_t$, which can be solved by means of the Newton method: given $\bm\mu\in\mathcal{P}$ and an initial guess $\mathbf{u}_h^{n,(0)}(\bm\mu)$, for $k\geq0$, find $\mathbf{\delta u}_h^{n,(k)}(\bm\mu)\in\mathbb{R}^{N_h}$ such that
\begin{equation}\label{eq:motionFOM} \left\{ \begin{array}{llll}
\mathbf{J}(\mathbf{u}_h^{n,(k)}(\bm\mu),t^{n};\bm\mu)\mathbf{\delta u}_h^{n,(k)}(\bm\mu) = - \mathbf{R}(\mathbf{u}_h^{n,(k)}(\bm\mu),t^{n};\bm\mu) \\
\mathbf{u}_h^{n,(k+1)}(\bm\mu) = \mathbf{u}_h^{n,(k)}(\bm\mu) + \mathbf{\delta u}_h^{n,(k)}(\bm\mu)
\end{array} \right.
\end{equation}
until suitable stopping criteria are fulfilled. Here, $\mathbf{u}_h^{n,(k)}(\bm\mu)$ represents the solution vector for a fixed parameter $\bm\mu$ computed at time step $t^{n}$ and Newton iteration $k$, while $\mathbf{R}\in\mathbb{R}^{N_h}$ and $\mathbf{J}\in\mathbb{R}^{N_h\times N_h}$ denote the residual vector and the corresponding Jacobian matrix, respectively. We refer to (\ref{eq:motionFOM}) as the high-fidelity, full-order model (FOM) for problem (\ref{eq:motion}). In particular, we rely on a Galerkin-finite element method (FEM) for space approximation, and  consider implicit finite difference schemes for time discretization, i.e.,
\begin{align*}
\partial_t\mathbf{u}_h(t^{n}) \approx \frac{\mathbf{u}_h^{n}-\mathbf{u}_h^{n-1}}{\Delta t}, &&
\partial_t^2\mathbf{u}_h(t^{n}) \approx \frac{\mathbf{u}_h^{n}-2\mathbf{u}_h^{n-1}+\mathbf{u}_h^{n-2}}{\Delta t^2},
\end{align*}
which do not require restrictions on $\Delta t$ \cite{quarteroni2013numerical}. 
The high-fidelity dimension $N_h$ is determined by the underlying mesh and the chosen FE polynomial order and can be extremely big whenever high accuracy is required for the problem at hand. 

To reduce the FOM numerical complexity, we introduce a projection-based reduced-order model (ROM), by relying on the reduced basis (RB) method  \cite{quarteroni2016reduced}. The idea of the RB method is to suitably select $N\ll N_h$ vectors of $\mathbb{R}^{N_h}$, forming the so-called RB matrix $\mathbf{V}\in\mathbb{R}^{N_h\times N}$, and to generate a reduced problem by performing a Galerkin projection of the FOM onto the subspace $V_N=\text{Col}(\mathbf{V})\subset\mathbb{R}^{N_h}$ generated by these vectors. This method relies on the assumption that the reduced-order approximation can be expressed as a linear combination of few, problem-dependent, basis functions, that is
\begin{equation*}
\mathbf{Vu}_N^{n}(\bm\mu) \approx \mathbf{u}_h^{n}(\bm\mu),
\end{equation*}
for $n=\,\dots,N_t$, where $\mathbf{u}_N^{n}(\bm\mu)\in\mathbb{R}^N$ denotes the vector of the ROM degrees of freedom at time $t^n\geq0$. The latter is obtained by imposing that the projection of the FOM residual computed on the ROM solution is orthogonal to the trial subspace (in the case of a Galerkin projection): given $\bm\mu\in\mathcal{P}$, at each time $t^{n}$, for $n=1,\dots,N_t$, we seek $\mathbf{u}_N^{n}(\bm\mu)\in\mathbb{R}^{N}$ such that 
\begin{equation*}
\mathbf{V}^T\mathbf{R}(\mathbf{Vu}_N^{n}(\bm\mu),t^{n};\bm\mu) = \mathbf{0}.
\end{equation*}
From now on, we will denote the reduced residual $\mathbf{V}^T\mathbf{R}$ and the corresponding Jacobian $\mathbf{V}^T\mathbf{JV}$ as $\mathbf{R}_N$ and $\mathbf{J}_N$, respectively. Then, the associated reduced Newton problem at time $t^{n}$ reads: given $\mathbf{u}_N^{n,(0)}(\bm\mu)$, for $k\geq0$, find $\mathbf{\delta u}_N^{n,(k)}(\bm\mu)\in\mathbb{R}^{N}$ such that
\begin{equation}\label{eq:motionROM} \left\{ \begin{array}{llll}
\mathbf{J}_N(\mathbf{V}\mathbf{u}_N^{n,(k)}(\bm\mu),t^{n};\bm\mu)\mathbf{\delta u}_N^{n,(k)}(\bm\mu) = - \mathbf{R}_N(\mathbf{V}\mathbf{u}_N^{n,(k)}(\bm\mu),t^{n};\bm\mu), \\
\mathbf{u}_N^{n,(k+1)}(\bm\mu) = \mathbf{u}_N^{n,(k)}(\bm\mu) + \mathbf{\delta u}_N^{n,(k)}(\bm\mu),
\end{array} \right.
\end{equation}
until a suitable stopping criterion is fulfilled. 


\subsection{Solution-space reduction: proper orthogonal decomposition}

In this section we provide an overview of the proper orthogonal decomposition (POD) technique used to compute the reduced basis $\mathbf{V}$ through the so-called method of snapshots \cite{POD,benner2015survey}. 
Let
\begin{equation*}
\mathcal{M}_{u_h} = \{\mathbf{u}_h^n(\bm\mu)\in\mathbb{R}^{N_h}~\lvert ~\bm\mu\in\mathcal{P}, ~n=1,\dots,N_t \}
\end{equation*}
be the (discrete) solution manifold identified by the image of $\mathbf{u}_h$, that is, the set of all the PDE solutions for $\bm\mu$ varying in the parameter space and $t^n$ in the partition of the time interval. Our goal is to approximate $\mathcal{M}_{u_h}$ with a reduced linear manifold, the {\itshape trial manifold}
\begin{equation*}
\mathcal{M}_{u_N}^{lin} = \{\mathbf{Vu}_N^n(\bm\mu)~\lvert ~\mathbf{u}_N^n(\bm\mu)\in\mathbb{R}^N, ~\bm\mu\in\mathcal{P},~n=1,\dots,N_t \}.
\end{equation*}
To do this, given $n_s<N_h$ sampled instances of  $\bm\mu\in\mathcal{P}$, we define the snapshots matrix 
\begin{equation*}
\mathbf{S} = \left[ \mathbf{u}_h(t^1;\bm\mu_1)~|~\dots~|~\mathbf{u}_h(t^{N_t};\bm\mu_1)~|~\dots~|~\mathbf{u}_h(t^1;\bm\mu_{n_s})~|~\dots~|~\mathbf{u}_h(t^{N_t};\bm\mu_{n_s}) \right]\in\mathbb{R}^{N_h\times n_s}
\end{equation*}
which contains the FOM solutions $\mathbf{u}_h(t^n;\bm\mu_\ell)$ as its columns. Sampling can be performed, e.g., through a latin hypercube sampling design, as well as through suitable low-discrepancy points sets. 

The POD basis $\mathbf{V}\in\mathbb{R}^{N_h\times N}$ spanning the subspace $V_N$ is obtained by performing the singular value decomposition of $\mathbf{S}$, 
\begin{equation*}
\mathbf{S}=\mathbf{U\Sigma Z}^T,
\end{equation*}
and then collecting the first $N$ columns of $\mathbf{U}$, corresponding to the $N$ largest singular values stored in the diagonal matrix $\bm\Sigma = \text{diag}(\sigma_1,\dots,\sigma_{r})\in\mathbb{R}^{n_s\times n_s}$, with $\sigma_1\geq\dots\geq\sigma_{r}\geq0$ and $r\leq N_h\land n_s$ being the rank of $\mathbf{S}$. The columns of the matrices $\mathbf{U}\in\mathbb{R}^{N_h\times n_s}$ and $\mathbf{Z}\in\mathbb{R}^{n_s\times N_h}$ correspond to the left and the right singular vectors of $\mathbf{S}$, respectively. This yields an orthonormal basis that, among all $N$-dimensional orthonormal basis $\mathbf{W}\in\mathcal{V}_N$, minimizes the least square error of snapshot reconstruction
\begin{align*}
&\lVert \mathbf{S} - \mathbf{V}\mathbf{V}^T\mathbf{S} \rVert_2^2 = \underset{\mathbf{W}\in\mathcal{V}_N}{\min} \lVert \mathbf{S} - \mathbf{W}\mathbf{W}^T\mathbf{S} \rVert_2^2 = \sigma_{N+1}^2, \\
&\lVert \mathbf{S} - \mathbf{V}\mathbf{V}^T\mathbf{S} \rVert_F^2 = \underset{\mathbf{W}\in\mathcal{V}_N}{\min} \lVert \mathbf{S} - \mathbf{W}\mathbf{W}^T\mathbf{S} \rVert_F^2 = \sum_{i=N+1}^r \sigma_i^2,
\end{align*}
where $\lVert \cdot \rVert_2$ and $\lVert \cdot \rVert_F$ are the Euclidean norm and the Frobenius norm, respectively. Hence, singular values' decay directly impacts on the size $N$, usually computed as the minimum integer satisfying  \vspace{-0.1cm}
\begin{equation}\label{eq:POD-dim}
RIC(N) = \frac{\sum_{\ell=1}^{N}\sigma_i^2}{\sum_{\ell=1}^{r}\sigma_i^2} \geq 1-\varepsilon_{POD}^2 \vspace{-0.1cm}
\end{equation}
for a given tolerance $\varepsilon_{POD}>0$. The POD technique constructs a low-dimensional optimal subspace of $\mathbb{R}^{N_h}$ retaining as much as possible of the snapshots relative information  content ($RIC$), provided that a sufficiently rich set of snapshots has been chosen. We summarize the POD technique in Algorithm \ref{alg:POD}.

\begin{algorithm}\caption{Proper orthogonal decomposition (POD)}\label{alg:POD}
	INPUT: $\mathbf{S}\in\mathbb{R}^{N_h\times n_s}$\\
	OUTPUT: $\mathbf{V}\in\mathbb{R}^{N_h\times N}$\\ \vspace{-1.1em}
	\begin{algorithmic}[1]
		\STATE Perform SVD of $\mathbf{S}$, i.e.,  $\mathbf{S=U\Sigma Z}^T$
		\STATE Select basis dimension $N$ as the minimum integer fulfilling condition (\ref{eq:POD-dim})
		\STATE Construct $\mathbf{V}$ collecting the first $N$ columns of $\mathbf{U}$
	\end{algorithmic}
\end{algorithm}

For the sake of efficiency, in our work we rely on a non-deterministic version of POD, exploiting the so-called randomized-SVD, see Algorithm \ref{alg:randSVD}. Randomization offers, in fact, a powerful tool for performing low-rank matrix approximation, especially when dealing with massive datasets. The randomized approach usually beats its classical competitors in terms of computational speed-up, accuracy and robustness \cite{halko2011finding}. 
The key idea of randomized SVD is to split the task of computing an approximated singular value decomposition of a given matrix into a first random stage, and a second deterministic one. The former exploits random sampling to construct a low-dimensional subspace that captures most of the action of the input matrix; the latter is meant to restrict the given matrix to this subspace and then manipulate the associated reduced matrix with classical deterministic algorithms, to obtain the desired low-rank approximations. This randomized approach is convenient when the snapshots matrix is high-dimensional, i.e. when $N_h$ and $n_s$ are large. In fact, finding the first $k$ dominant singular-values for a dense input matrix of dimension $N_h\times n_s$, requires $O(N_hn_s\log(k))$ floating-point operations for a randomized algorithm, in contrast with $O(N_hn_sk)$ flops for a classical one.  	 

	\begin{algorithm}\caption{Randomized-SVD}\label{alg:randSVD}
		INPUT: $\mathbf{S\in\mathbb{R}}^{m\times n}$, target rank $k\in\mathbb{N}$\\
		OUTPUT: $\mathbf{U\Sigma Z}^T\approx\mathbf{S}$\\ \textbf{\textit{stage 1}}
		\begin{algorithmic}[1]
			\STATE Generate a Gaussian matrix $\bm\Theta\in\mathbb{R}^{N_h\times k}$
			\STATE Compute $\mathbf{Q}\in\mathbb{R}^{N_h\times k}$ whose columns form an orthonormal basis for the range of $\mathbf{S}\bm\Theta$ and such that 
			\begin{equation*}
			\lVert \mathbf{S} - \mathbf{QQ}^T\mathbf{S}\rVert_2 \leq \min_{rank(\mathbf{X})\leq k} \lVert \mathbf{S} - \mathbf{X}\rVert_2,
			\end{equation*}
			e.g., using the QR factorization.
		\end{algorithmic}
		\textbf{\textit{stage 2}}
		\begin{algorithmic}[1]
			\STATE Form $\widetilde{\mathbf{S}}=\mathbf{Q}^T\mathbf{S}\in\mathbb{R}^{k\times n}$
			\STATE Compute SVD of $\widetilde{\mathbf{S}}=\widetilde{\mathbf{U}}\bm\Sigma \mathbf{Z}^T$
			\STATE Set $\mathbf{U}=\mathbf{Q}\widetilde{\mathbf{U}}$
		\end{algorithmic}
	\end{algorithm}		
	
	\begin{remark}
	Note that Algorithm~\ref{alg:randSVD} can be also adapted to solve the following problem: given a target error tolerance $\varepsilon>0$, find $k=k(\varepsilon)$ and $\mathbf{Q}\in\mathbb{R}^{N_h\times k}$ satisfying 
	$\lVert \mathbf{S} - \mathbf{QQ}^T\mathbf{S}\rVert_2 \leq \varepsilon$. 
 \end{remark}


\subsection{Hyper-reduction: the discrete empirical interpolation method}\label{sec:DEIM}
 
In the case of parametrized PDEs featuring nonaffine dependence on the parameter and/or nonlinear (high-order polynomial or nonpolynomial) dependence on the field variable, a further level of reduction, known as {\itshape hyper-reduction}, must be introduced \cite{grepl2007efficient,negri2015efficient}. Note that if nonlinearities only include quadratic (or, at most, cubic) terms and do not feature any parameter dependence, assembling of nonlinear terms in the ROM can be performed by projection of the corresponding FOM quantities, once and for all \cite{gobat2021reduced}.

For the case at hand, the residual $\mathbf{R}_N$ and the Jacobian $\mathbf{J}_N$ appearing in the reduced Newton system (\ref{eq:motionROM}) depend on the solution at the previous iteration and, therefore, must be computed at each step $k\geq0$. It follows that, for any new instance of the parameter $\bm\mu$, we need to assemble the high-dimensional FOM-arrays before projecting them onto the reduced subspace, entailing a computational complexity which is still of order $N_h$. To setup an efficient offline-online computational splitting, an approximation of the nonlinear operators that is independent of the high-fidelity dimension is required.  

Several techniques have been employed to provide this further level of approximation \cite{barrault2004empirical, chaturantabut2010nonlinear, astrid2008missing, carlberg2011efficient, farhat2015structure}; among these, DEIM has been successfully applied to stationary or quasi-static nonlinear mechanical problems \cite{bonomi2017reduced,ghavamian2017pod}. Its key idea is to replace the nonlinear arrays in (\ref{eq:motionROM}) with a collateral reduced basis expansion, computed through an inexpensive interpolation procedure. In this framework, the high-dimensional residual $\mathbf{R}(\bm\mu)$ is projected onto a reduced subspace of dimension $m<N_h$ spanned by a basis $\bm\Phi_\mathcal{R}\in\mathbb{R}^{N_h\times m}$ \vspace{-0.05cm}
\begin{equation*}
\mathbf{R}(\bm\mu) \approx \bm\Phi_\mathcal{R}\mathbf{r}(\bm\mu), \vspace{-0.05cm}
\end{equation*}
where $\mathbf{r}(\bm\mu)\in\mathbb{R}^m$ is the vector of the unknown amplitudes. The matrix $\bm\Phi_\mathcal{R}$ can be precomputed offline by performing POD on a set of high-fidelity residuals collected when solving (\ref{eq:motionROM}) for $n_s'$ training input parameters \vspace{-0.05cm}
\begin{equation*}
\mathbf{S}_{\bm\rho} =\left[\mathbf{R}(\mathbf{Vu}_N^{n,(k)}(\bm\mu_\ell),t^{n};\bm\mu_\ell)), k\geq0\right]_{n=1,\dots,N_t}^{\ell=1,\dots,n_s'}. \vspace{-0.05cm}
\end{equation*}
The unknown parameter-dependent coefficient $\mathbf{r}(\bm\mu)$ is obtained online by collocating the approximation at the $m$ components selected by a greedy procedure, that is \vspace{-0.05cm}
\begin{equation*}
\mathbf{P}^T\mathbf{R}(\bm\mu) \approx \mathbf{P}^T\bm\Phi_\mathcal{R}\mathbf{r}(\bm\mu) \implies \mathbf{r}(\bm\mu) = (\mathbf{P}^T\bm\Phi_\mathcal{R})^{-1}\mathbf{P}^T\mathbf{R}(\bm\mu), \vspace{-0.05cm}
\end{equation*}
where $\mathbf{P}\in\mathbb{R}^{N_h\times m}$ is the boolean matrix associated with the interpolation constraints. We thus define the hyper-reduced residual vector as \vspace{-0.05cm}
\begin{equation*}
\mathbf{R}_{N,m}(\bm\mu) := \mathbf{V}^T\bm\Phi_\mathcal{R}(\mathbf{P}^T\bm\Phi_\mathcal{R})^{-1}\mathbf{P}^T\mathbf{R}(\bm\mu) \approx \mathbf{V}^T\mathbf{R}(\bm\mu). \vspace{-0.05cm}
\end{equation*}
To avoid confusion, we recall that $\mathbf{R}_N = \mathbf{V}^T\mathbf{R}$, so that $\mathbf{R}_{N,m}\approx\mathbf{R}_N$. Finally, the associated Jacobian approximation $\mathbf{J}_{N,m}(\bm\mu)$ can be computed as the derivative of $\mathbf{R}_{N,m}(\bm\mu)$ with respect to the reduced displacement, obtaining \vspace{-0.05cm}
\begin{align*}
\mathbf{J}_{N,m}(\bm\mu) = \mathbf{V}^T{\bm\Phi_\mathcal{R}}(\mathbf{P}^T{\bm\Phi_\mathcal{R}})^{-1}	 \mathbf{P}^T\mathbf{J}(\bm\mu)\mathbf{V}, \vspace{-0.05cm}
\end{align*}
or by relying on the so-called matrix DEIM (MDEIM) algorithm \cite{negri2015efficient}, as done in \cite{bonomi2017reduced,manzoni2018reduced}.

\begin{figure}[t!]
	\centering
	\includegraphics[scale=0.475]{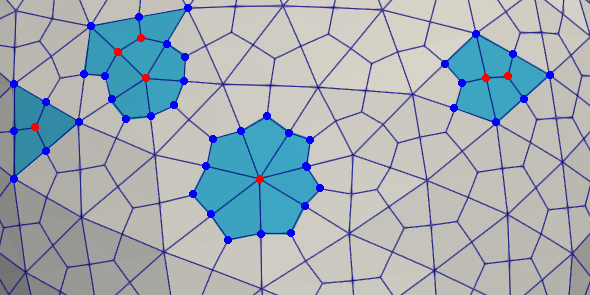}
	\caption{Sketch of a reduced mesh for an hexahedral computational grid in a two-dimensional case. The red dots represent the points selected by the DEIM algorithm and, together with the blue ones, represent the vertices of the elements (light blue) of the reduced mesh.}
	\label{fig:redMesh}
\end{figure}

However, the application of DEIM in this setting can be rather inefficient, especially when turning to complex problem which require a high number of residual basis, thus interpolation points, to ensure the convergence of the hyper-reduced Newton system
\begin{equation*}
\left\{ \begin{array}{llll}
\mathbf{J}_{N,m}(\mathbf{V}\mathbf{u}_N^{n,(k)}(\bm\mu),t^{n};\bm\mu)\mathbf{\delta u}_N^{n,(k)}(\bm\mu) = - \mathbf{R}_{N,m}(\mathbf{V}\mathbf{u}_N^{n,(k)}(\bm\mu),t^{n};\bm\mu) \\
\mathbf{u}_N^{n,(k+1)}(\bm\mu) = \mathbf{u}_N^{n,(k)}(\bm\mu) + \mathbf{\delta u}_N^{n,(k)}(\bm\mu).
\end{array} \right. 
\end{equation*}
In fact, the $m$ points selected by the DEIM algorithm correspond to a subset of nodes of the computational mesh, which, together with the neighboring nodes (i.e. those sharing the same cell), form the so-called \textit{reduced mesh}, see, e.g., the sketch reported in Figure~\ref{fig:redMesh}. Since the entries of any FE-vector are associated with the degrees of freedom (dofs) of the problem, $\mathbf{P}^T\mathbf{R}(\bm\mu)$ is computed by integrating the residual only on the quadrature points belonging to the reduced mesh, which, nevertheless, can be rather dense.

A modification of the DEIM algorithm, the so-called unassembled DEIM (UDEIM), has been proposed in \cite{tiso2013discrete} to exploit the sparsity of the problem and minimize the number of element function calls. However, a high number of nonlinear function evaluations is still required when the number of magic points is sufficiently big. Indeed,  DEIM-based affine approximations are effective, in terms of computational costs, provided that few entries of the nonlinear terms can be cheaply computed; however, this situation does not occur neither for dynamical systems arising from the linearization of a nonlinear system around a steady state, nor when dealing with global nonpolynomial nonlinearities. 

In this paper, we propose an alternative technique to perform hyper-reduction, which is independent of the underlying mesh and relies on a deep neural network architecture to approximate reduced residual vectors and Jacobian matrices. The introduction of a surrogate model to perform operator approximation is justified by the fact that, often, most of the CPU time needed online for each new parameter instance is required by DEIM for assembling arrays such as residual vectors or corresponding Jacobian matrices on the reduced mesh.


\section{Operator approximation: a deep learning-based technique (Deep-HyROMnet)}\label{sec:Deep-HyROMnet}

To recover the offline-online efficiency of the RB method, overcoming the need to assemble the nonlinear arrays onto the computational mesh as in the case of the DEIM, we present a novel projection-based method which relies on DNNs for the approximation of the nonlinear terms. We refer to this strategy as to a \textit{hyper-reduced order model enhanced by deep neural networks} (Deep-HyROMnet). Our goal is the efficient numerical approximation of the whole sets
\begin{align*}
\mathcal{M}_{R_N} &= \{\mathbf{R}_N(\mathbf{Vu}_N^{n,(k)}(\bm\mu),t^n;\bm\mu)\in\mathbb{R}^{N} ~\lvert ~n=1,\dots,N_t, ~k\geq0, ~\bm\mu\in\mathcal{P} \}, \medskip \\
\mathcal{M}_{J_N} &= \{\mathbf{J}_N(\mathbf{Vu}_N^{n,(k)}(\bm\mu),t^n;\bm\mu)\in\mathbb{R}^{N\times N} ~\lvert ~n=1,\dots,N_t, ~k\geq0, ~\bm\mu\in\mathcal{P} \},
\end{align*}
which we refer to as the \textit{reduced residual manifold} and \textit{reduced Jacobian manifold}, respectively, in a way that depends only on the ROM dimension $N$ and on the number of parameters $P$. To achieve this task, we employ the DNN architecture developed in \cite{fresca2021comprehensive} for the DL-ROM techniques. It is worthy to note that, except for the approximation error of the reduced nonlinear operators, the proposed Deep-HyROMnet approach is a physics-based method and that the computed solution satisfies the nonlinear equation of the problem under investigation, up to a further approximation of ROM residual and Jacobian arrays -- thus, similarly to what happened for a POD-Galerkin-DEIM ROM. The main idea of the deep learning-based operator approximation approach that replaces the DEIM in our new Deep-HyROMnet strategy, is to learn the following input-to-residual and input-to-Jacobian maps, respectively:
\begin{align*}
\bm\rho_N\colon&(\bm\mu,t^n,k)\longmapsto {\bm\rho}_N(\bm\mu,t^n,k) \approx \mathbf{R}_N(\mathbf{Vu}_N^{n,(k)}(\bm\mu),t^n;\bm\mu),\\
\bm\iota_N\colon&(\bm\mu,t^n,k)\longmapsto {\bm\iota}_N(\bm\mu,t^n,k) \approx \mathbf{J}_N(\mathbf{Vu}_N^{n,(k)}(\bm\mu),t^n;\bm\mu),
\end{align*}
provided $(\bm\mu,t^n,k)\in\mathcal{P}\times\{t^1,\dots,t^{N_t}\}\times\mathbb{N}^+$, and to finally replace the linear system in (\ref{eq:motionROM}) with
\begin{equation*}
{\bm\iota}_N(\bm\mu,t^n,k)\mathbf{\delta u}^{n,(k)}(\bm\mu) = - {\bm\rho}_N(\bm\mu,t^n,k).
\end{equation*}

Hence, Deep-HyROMnet aims at approximating the residual vector and the Jacobian matrix obtained after their projection onto the reduced space of dimension $N\ll N_h$. Indeed, performing POD-Galerkin on the solution space allows to severely reduce the problem dimension from $N_h$ to $N$ and, hence, to ease the learning task with respect the reconstruction of the full-order $\mathbf{R}$ and $\mathbf{J}$.

\begin{remark}
As an alternative to Newton iterative scheme, we can rely on Broyden's method \cite{Broyden}, which belongs to the class of quasi-Newton methods. This allows to avoid the computation of the Jacobian matrix at each iteration $k\geq0$ by relying on rank-one updates, based on residuals computed at previous iterations. However, we are able to compute  
Jacobian matrices  
 very efficiently using automatic differentiation (AD), so that the computational burden is the assembling of residual vectors. For this reason, in this paper we will focus on the Newton method only, that is, the solution of problem (\ref{eq:motionROM}).
\end{remark}

To summarize, in the case of the Newton approach, we end up with the following reduced problem: given $\bm\mu\in\mathcal{P}$ and, for $n=1,\dots,N_t$, the initial guess $\mathbf{u}_N^{n,(0)}(\bm\mu) = \mathbf{u}_N^{n-1}(\bm\mu)$, find $\delta\mathbf{u}_N^{n,(k)}\in\mathbb{R}^N$ such that, for $k\geq0$, 
\begin{equation}\label{eq:motionDeep} \left\{ \begin{array}{llll}
\bm\iota_N(\bm\mu,t^{n},k) \delta\mathbf{u}_N^{n,(k)}(\bm\mu) = - \bm\rho_N(\bm\mu,t^{n},k), \\
\mathbf{u}_N^{n,(k+1)}(\bm\mu) = \mathbf{u}_N^{n,(k)}(\bm\mu) + \delta\mathbf{u}_N^{n,(k)}(\bm\mu),
\end{array} \right.
\end{equation}
until $\lVert\bm\rho_N(\bm\mu,t^n,k)\rVert_2 / \lVert\bm\rho_N(\bm\mu,t^n,0)\rVert_2 < \varepsilon$, where $\varepsilon>0$ is a given tolerance. In Algorithms~\ref{alg:Deep-HyROMnet_offline} and \ref{alg:Deep-HyROMnet_online}, we report a summary of the offline and online stages of Deep-HyROMnet, respectively.
\begin{algorithm}[H]
	\caption{Deep-HyROMnet for nonlinear time-dependent problems. Offline stage.}
	\label{alg:Deep-HyROMnet_offline}
	INPUT: $\bm\mu_\ell$, for $\ell=1,\dots,n_s$, and $\bm\mu_{\ell'}$, for $\ell'=1,\dots,n_s'$\\
	OUTPUT: $\mathbf{V}\in\mathbb{R}^{N_h\times N}$
	\begin{algorithmic}[1]
		\STATE \textbf{for} $\ell=1,\dots,n_s$ \textbf{do}
		\STATE $\quad$ \textbf{for} $n=1,\dots,N_t$ \textbf{do}
		\STATE $\quad\quad$ \textbf{for} $k\geq0$ \textbf{until convergence do}
		\STATE $\quad\quad\quad$ Assemble and solve problem (\ref{eq:motionFOM})
		\STATE $\quad\quad\quad$ Collect $\mathbf{S}_u \leftarrow\mathbf{S}_u\cup \left[\mathbf{u}_h^{n,(k)}(\bm\mu_\ell)\right]$ column-wise
		\STATE Construct $\mathbf{V}=\text{POD}(\mathbf{S}_u, \varepsilon_{POD})$ (see Algorithm~\ref{alg:POD})
		\STATE \textbf{for} $\ell'=1,\dots,n_s'$\textbf{ do}
		\STATE $\quad$ \textbf{for} $n=1,\dots,N_t$ \textbf{do}
		\STATE $\quad\quad$ \textbf{for} $k\geq0$ \textbf{until convergence do}	
		\STATE $\quad\quad\quad$ Assemble and solve reduced problem (\ref{eq:motionROM})
		\STATE $\quad\quad\quad$ Collect $\mathbf{S}_{\bm\rho} \leftarrow\mathbf{S}_{\bm\rho}\cup \left[\mathbf{R}_N(\mathbf{V}\mathbf{u}_N^{n,(k)}(\bm\mu_{\ell'}),t^{n};\bm\mu_{\ell'})\right]$ column-wise
		\STATE $\quad\quad\quad$ Collect $\mathbf{S}_{\bm\iota} \leftarrow\mathbf{S}_{\bm\iota}\cup \left[\mathbf{J}_N(\mathbf{V}\mathbf{u}_N^{n,(k)}(\bm\mu_{\ell'}),t^{n};\bm\mu_{\ell'})\right]$ column-wise
		\STATE Train the DNNs (see Algorithm~\ref{alg:DLROM1})
	\end{algorithmic}
\end{algorithm}

\begin{algorithm}[H]
	\caption{Deep-HyROMnet for nonlinear time-dependent problems. Online stage.}
	\label{alg:Deep-HyROMnet_online}	
	INPUT: $\bm\mu\in\mathcal{P}$\\
	OUTPUT: $\mathbf{Vu}_N^n(\bm\mu)\in\mathbb{R}^{N_h}$, for $n=1,\dots,N_t$\\ \vspace{-1.1em}
	\begin{algorithmic}[1]
		\STATE \textbf{for} $n=0,\dots,N_t-1$ \textbf{do}
		\STATE $\quad$ \textbf{for }$k\geq0$ \textbf{until convergence do}
		\STATE $\quad\quad$ Compute $\bm\rho_N(\bm\mu,t^n,k)$ and $\bm\iota_N(\bm\mu,t^n,k)$ (see Algorithm~\ref{alg:DLROM2})
		\STATE $\quad\quad$ Solve hyper-reduced problem (\ref{eq:motionDeep})
		\STATE Recover $\mathbf{Vu}_N^n(\bm\mu)$, for $n=1,\dots,N_t$
	\end{algorithmic}
\end{algorithm}


\subsection{DL-ROM-based neural network}\label{sec:DLROM}

For the sake of generality, we will focus on the DNN-based approximation of the reduced residual vector only, that is 
\begin{equation*}
{\bm\rho}_N(\bm\mu,t^n,k) \approx\mathbf{R}_N(\mathbf{V}\mathbf{u}_N^{n,(k)}(\bm\mu),t^n;\bm\mu)\in\mathbb{R}^N.
\end{equation*}
In fact, by relying on a suitable transformation, we can easily write the Jacobian matrix as a vector of dimension $N^2$ and apply the same procedure described in the following for the residual vector also in the case of the Jacobian matrix. In particular, we define the transformation
\begin{equation*}
vec\colon\mathbb{R}^{N\times N}\rightarrow\mathbb{R}^{N^2}, \quad vec(\mathbf{J}_N(\mathbf{V}\mathbf{u}_N^{n,(k)}(\bm\mu),t^n;\bm\mu)) = \mathbf{j}_N(\mathbf{V}\mathbf{u}_N^{n,(k)}(\bm\mu),t^n;\bm\mu),
\end{equation*}
which consists in stacking the columns of $\mathbf{J}_N(\mathbf{V}\mathbf{u}_N^{n,(k)}(\bm\mu),t^n;\bm\mu)$ in a vector on which is then applied the DL-ROM technique, thus obtaining 
\begin{equation*}
\widetilde{\bm\iota}_N(\bm\mu,t^n,k)\approx \mathbf{j}_N(\mathbf{V}\mathbf{u}_N^{n,(k)}(\bm\mu),t^n;\bm\mu)\in\mathbb{R}^{N^2}.
\end{equation*}
Finally, we revert the $vec$ operation, so that ${\bm\iota}_N(\bm\mu,t^n,k) = vec^{-1}(\widetilde{\bm\iota}_N(\bm\mu,t^n,k))$.

We thus aim at efficiently approximating the whole set $\mathcal{M}_{R_N}$
by means of the reduced residual trial manifold, defined as
\begin{equation*}
\mathcal{M}_{\rho_N} = \{ {\bm\rho}_N(\bm\mu,t^n,k) ~|~\bm\mu\in\mathcal{P}, ~n=1,\dots,N_t, ~k\geq0\}\subset\mathbb{R}^{N}.
\end{equation*}
The DL-ROM approximation of the ROM residual $\mathbf{R}_N(\mathbf{V}\mathbf{u}_N^{n,(k)}(\bm\mu),t^n;\bm\mu)$ takes the form
\begin{equation*}
\bm\rho_N(\bm\mu,t^n,k) = \widetilde{\mathbf{R}}_N(\bm\mu,t^n,k;\bm\theta_{DF},\bm\theta_{D}) = \mathbf{f}^D_N(\bm\phi_q^{DF}(\bm\mu,t^n,k;\bm\theta_{DF});\bm\theta_{D})
\end{equation*}
where
\begin{itemize}
	\item $\bm\phi_q^{DF}(\cdot~;\bm\theta_{DF})\colon\mathbb{R}^{P+2}\rightarrow\mathbb{R}^q$ such that 
	\begin{equation*}
	\bm\phi_q^{DF}(\bm\mu,t^n,k;\bm\theta_{DF}) = \mathbf{R}_q(\bm\mu,t^n,k;\bm\theta_{DF})
	\end{equation*} is a deep feedforward neural network (DFNN), consisting in the subsequent composition of a nonlinear activation function, applied to a linear transformation of the input, multiple times. Here, $\bm\theta_{DF}$ denotes the vector of parameters of the DFNN, collecting all the corresponding weights and biases of each layer and $q$ is as close as possible to the input size $P+2$;
	\item $\mathbf{f}^D_N(\cdot~;\bm\theta_{D})\colon\mathbb{R}^q\rightarrow\mathbb{R}^N$ such that
	\begin{equation*}
	\mathbf{f}^D_N(\mathbf{R}_q(\bm\mu,t^n,k;\bm\theta_{DF});\bm\theta_{D}) = \widetilde{\mathbf{R}}_N(\bm\mu,t^n,k;\bm\theta_{DF},\bm\theta_{D})
	\end{equation*}
	is the decoder function of a convolutional autoencoder (CAE), obtained as the composition of several layers (some of which are convolutional), depending upon a vector $\bm\theta_{D}$ collecting all the corresponding weights and biases.
\end{itemize}
The encoder function of the CAE is exploited, during the training stage only, to map the reduced residual $\mathbf{R}_N(\mathbf{V}\mathbf{u}_N^{n,(k)}(\bm\mu),t^n;\bm\mu)$ associated to $(\bm\mu,t^n,k)$ onto a low-dimensional representation
\begin{equation*}
\mathbf{f}^E_q(\mathbf{R}_N(\mathbf{V}\mathbf{u}_N^{n,(k)}(\bm\mu),t^n;\bm\mu);\bm\theta_{E}) = \widetilde{\mathbf{R}}_q(\bm\mu,t^n,k;\bm\theta_{E}),
\end{equation*}
where $\mathbf{f}^E_q(\cdot~;\bm\theta_{E})\colon\mathbb{R}^N\rightarrow\mathbb{R}^q$ denotes the encoder function, depending upon a vector $\bm\theta_{E}$ of parameters. 

\begin{remark}\label{rmk:zero-padded}
	We point out that the input of the encoder function, that is, the reduced residual vector $ \mathbf{R}_N$, is reshaped into a square matrix by rewriting its elements in row-major order, thus obtaining $\mathbf{R}_N^{reshape}\in\mathbb{R}^{\sqrt{N}\times \sqrt{N}}$. If $N$ is not a square, the input $\mathbf{R}_N$ is zero-padded as explained in \cite{Goodfellow-et-al-2016}, and the additional elements are subsequently discarded.	
\end{remark}

\begin{figure}
	\centering
	\includegraphics[width=0.95\textwidth]{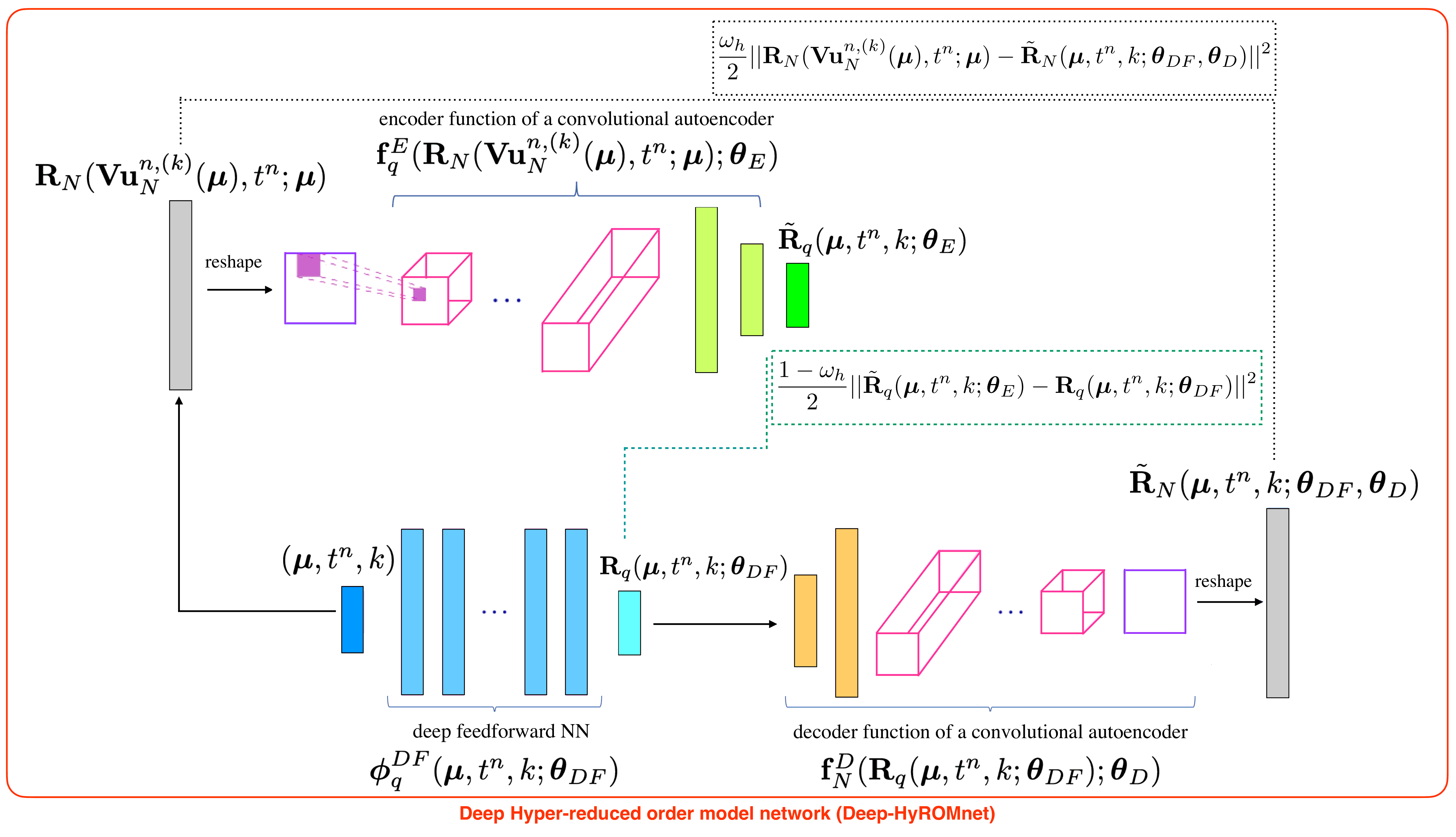}
	\caption{DNN architecture used during the training phase for the reduced residual vector.}
	\label{fig:DNN}
\end{figure}

Regarding the prediction of the reduced residual for new unseen instances of the inputs, given $\bm\mu_{test}\in\mathcal{P}$, computing the DL-ROM approximation of $\mathbf{R}_N(\mathbf{V}\mathbf{u}_N^{n,(k)}(\bm\mu),t^n;\bm\mu_{test})$, for any possible $n=1,\dots,N_t$ and $k\geq0$, corresponds to the testing stage of a DFNN and of the decoder function of a convolutional AE; thus, at testing time, we discard the encoder function. The architecture used during the training stage is reported in Figure~\ref{fig:DNN}, whereas, during the testing phase, the encoder function is discarded.

Let us define the reduced residual snapshots matrix $\mathbf{S}_{\bm\rho}\in\mathbb{R}^{N\times N_{train}}$, with $N_{train}=n_s'N_tN_k$, as
\begin{equation*}	
\mathbf{S}_{\bm\rho} = \left[\mathbf{R}_N(\mathbf{V}\mathbf{u}_N^{n,(k)}(\bm\mu),t^n;\bm\mu_\ell)\right]_{\ell=1,\dots,n_s',n=1,\dots,N_t,k\geq0},
\end{equation*}
that is, the matrix collecting column-wise ROM residuals computed for $n_s'$ sampled parameters $\bm\mu_\ell\in\mathcal{P}$, at different time instances $t^1,\dots,t^{N_t}$ and for each Newton iteration $k\geq0$, and the parameter matrix $\mathbf{M}\in\mathbb{R}^{(P+2)\times N_{train}}$ of the corresponding triples
\begin{equation*}
\mathbf{M} = \left[\left(\bm\mu_\ell,t^n,k\right)\right]_{\ell=1,\dots,n_s',n=1,\dots,N_t,k\geq0}.
\end{equation*}
The training stage consists in solving the following optimization problem in the weights variable $\bm\theta = (\bm\theta_{E},\bm\theta_{DF},\bm\theta_{D})$:
\begin{equation*}
 \mathcal{J}(\bm\theta) =  \dfrac{1}{N_{train}}\sum_{\ell=1}^{n_s'}\sum_{n=1}^{N_t}\sum_{k=0}^{N_k}\mathcal{L}(\bm\mu_\ell,t^n,k;\bm\theta) \rightarrow \underset{\bm\theta}{\min}
\end{equation*}
where
\begin{equation}\label{eq:loss}
\begin{aligned}
\mathcal{L}(\bm\mu_\ell,t^n,k;\bm\theta) = &\dfrac{\omega_h}{2} \lVert \mathbf{R}_N(\mathbf{V}\mathbf{u}_N^{n,(k)}(\bm\mu_\ell),t^n;\bm\mu_\ell) -  \widetilde{\mathbf{R}}_N(\bm\mu_\ell,t^n,k;\bm\theta_{DF},\bm\theta_{D}) \rVert^2 \\
& + \dfrac{1-\omega_h}{2} \lVert \widetilde{\mathbf{R}}_q(\bm\mu_\ell,t^n,k;\bm\theta_{E}) - \mathbf{R}_q(\bm\mu_\ell,t^n,k;\bm\theta_{DF}) \rVert^2,
\end{aligned}
\end{equation}
with $\omega_h\in[0,1]$. The loss function (\ref{eq:loss}) combines the reconstruction error, i.e. the error between the ROM residual and the DL-ROM approximation, and the error between the intrinsic coordinates and the output of the encoder. The training stage of the DNN involved in Deep-HyROMnet is detailed in Algorithm~\ref{alg:DLROM1}; in particular, we denote by $\alpha$ the training-validation splitting fraction, by $\eta$ the starting learning rate, by $N_b$ the batch size, by $n_b = (1-\alpha)N_{train}/N_b$ the number of minibatches and by $N_e$ the maximum number of epochs. The testing stage of the DNN is detailed in Algorithm~\ref{alg:DLROM2}. See, e.g., \cite{fresca2021comprehensive, fresca2021pod} for further details.

\begin{algorithm}
	\caption{Training stage for the DNN, based on Algorithm~1 of \cite{fresca2021comprehensive}}
	\label{alg:DLROM1}
	INPUT: $\mathbf{M}\in\mathbb{R}^{(P+2)\times N_{train}}$, $\mathbf{S}_{\bm\rho}\in\mathbb{R}^{N\times N_{train}}$, $\alpha$, $\eta$, $N_b$, $n_b$, $N_e$, early-stopping criterion \\
	OUTPUT: $\bm\theta^* = (\bm\theta_E^*,\bm\theta_{DF}^*,\bm\theta_D^*)$ (optimal)
	\begin{algorithmic}[1]
		\STATE Randomly shuffle $\mathbf{M}$ and $\mathbf{S}$
		\STATE Split data in $\mathbf{M} = \left[\mathbf{M}^{train},\mathbf{M}^{val}\right]$ and $\mathbf{S}_{\bm\rho} = \left[\mathbf{S}_{\bm\rho}^{train},\mathbf{S}_{\bm\rho}^{val}\right]$ (according to $\alpha$)
		\STATE Normalize $\mathbf{M}$ and $\mathbf{S}$ according to (\ref{eq:standardization})
		\STATE Randomly initialize $\bm\theta^0 = (\bm\theta_E^0,\bm\theta_{DF}^0,\bm\theta_D^0)$
		\STATE $n_e = 0$
		\STATE \textbf{while} $\lnot$early-stopping \textbf{and} $n_e\leq N_e$ \textbf{do}
		\STATE $\quad$ \textbf{for} $k=1,\dots,n_b$ \textbf{do}
		\STATE $\quad\quad$ Sample a minibatch $(\mathbf{M}^{batch},\mathbf{S}^{batch})\subset(\mathbf{M}^{train},\mathbf{S}^{train})$
		\STATE\label{state:start} $\quad\quad$ $\mathbf{S}^{batch} = reshape(\mathbf{S}^{batch})$
		\STATE $\quad\quad$ $\widetilde{\mathbf{S}}^{batch}_q(\bm\theta_E^{n_bn_e+k}) = \mathbf{f}^E_q(\mathbf{S}^{batch};\bm\theta_E^{n_bn_e+k})$
		\STATE $\quad\quad$ $\mathbf{S}^{batch}_q(\bm\theta_{DF}^{n_bn_e+k}) = \bm\phi^{DF}_q(\mathbf{M}^{batch};\bm\theta_{DF}^{n_bn_e+k})$ 
		\STATE $\quad\quad$ $\widetilde{\mathbf{S}}^{batch}_N(\bm\theta_{DF}^{n_bn_e+k},\bm\theta_D^{n_bn_e+k}) = \mathbf{f}^D_N(\mathbf{S}^{batch}_q(\bm\theta_{DF}^{n_bn_e+k});\bm\theta_D^{n_bn_e+k})$
		\STATE\label{state:end} $\quad\quad$ $\widetilde{\mathbf{S}}^{batch}_N = reshape(\widetilde{\mathbf{S}}^{batch}_N)$
		\STATE $\quad\quad$ Accumulate loss (\ref{eq:loss}) on $(\mathbf{M}^{batch},\mathbf{S}^{batch})$ and compute $\hat{\nabla}_{\bm\theta}\mathcal{J}$
		\STATE $\quad\quad$ $\bm\theta^{n_bn_e+k+1} = \text{ADAM}(\eta,\hat{\nabla}_{\bm\theta}\mathcal{J},\bm\theta^{n_bn_e+k})$
		\STATE $\quad$ Repeat instructions \ref{state:start}-\ref{state:end} on $(\mathbf{M}^{val},\mathbf{S}^{val})$ to evaluate early-stopping criterion
		\STATE $\quad$ $n_e = n_e+1$
	\end{algorithmic}
\end{algorithm}
\begin{algorithm}
	\caption{Testing stage for the DNN, based on Algorithm~2 of \cite{fresca2021comprehensive}}
	\label{alg:DLROM2}
	INPUT: $(\bm\mu,t^n,k)\in\mathcal{P}\times\{t^1,\dots,t^{N_t}\}\times\mathbb{N}^+$, $(\bm\theta_{DF}^*,\bm\theta_D^*)$ (optimal)\\
	OUTPUT: $\widetilde{\mathbf{S}}_N$ (i.e. $\bm\rho_N(\bm\mu,t^n,k)$ or $\bm\iota_N(\bm\mu,t^n,k)$)
	\begin{algorithmic}[1]
		\STATE $\mathbf{S}_q(\bm\theta_{DF}^*) = \bm\phi^{DF}_q(\bm\mu,t^n,k;\bm\theta_{DF}^*)$
		\STATE $\widetilde{\mathbf{S}}_N(\bm\theta_{DF}^*,\bm\theta_D^*) = \mathbf{f}^D_N(\mathbf{S}_q(\bm\theta_{DF}^*);\bm\theta_D^*)$
		\STATE $\widetilde{\mathbf{S}}_N = reshape(\widetilde{\mathbf{S}}_N)$
	\end{algorithmic}
\end{algorithm}
\begin{remark}
	Differently from the scaling techniques used in \cite{fresca2021comprehensive,fresca2021pod}, which are based on a min-max procedure, we standardize the input and output of the DNN as follows. After splitting the data into training and validation sets according to a user-defined training-validation splitting fraction, $\mathbf{M} = \left[\mathbf{M}^{train},\mathbf{M}^{val}\right]$ and $\mathbf{S}_{\bm\rho} = \left[\mathbf{S}_{\bm\rho}^{train},\mathbf{S}_{\bm\rho}^{val}\right]$, we define for each row of the training set the corresponding mean and standard deviation 
	\begin{equation*}
	M_{mean}^i = \dfrac{1}{N_{train}} \sum_{j=1}^{N_{train}} M_{ij}^{train}, \quad M_{sd}^i = \sqrt{\dfrac{1}{N_{train}-1} \sum_{j=1}^{N_{train}} (M_{ij}^{train}-M_{mean}^i)^2},
	\end{equation*}
	so that parameters are normalized by applying the following transformation
	\begin{equation}\label{eq:standardization}
	M_{ij}^{train} \mapsto \dfrac{M_{ij}^{train}-M_{mean}^i}{M_{sd}^i}, \quad i=1,\dots,P+2, \quad j=1,\dots,N_{train}
	\end{equation}
	that is, each feature of the training parameter matrix is standardized. The same procedure is applied to the training snapshots matrix $\mathbf{S}_{\bm\rho}^{train}$ by replacing $M^i_{*}$ with $S^i_{*}$, where $*\in\{mean,sd\}$ respectively. Transformation (\ref{eq:standardization}) is applied to the validation and testing sets as well, but considering the mean and the standard deviation computed over the training set. In order to rescale the reconstructed solution to the original values, we apply the inverse transformation.
\end{remark}


\section{Numerical results}\label{sec:tests}

In this Section, we investigate the performances of Deep-HyROMnet on different applications related to the parametrized nonlinear time-dependent PDE problems, focusing on structural mechanics. In particular, we consider {\em (i)} a series of structural tests on a rectangular beam, with different loading conditions and a simple nonlinear constitutive law, and then {\em (ii)} a test case on an idealized left ventricle geometry, simulating cardiac contraction. In the following subsection we formulate both these problems in the framework of nonlinear elastodynamics. 

\subsection{Nonlinear elastodynamics} \label{sec:nonlinear_elastodynamics}

Let us consider a continuum body $\mathcal{B}$ embedded in a three-dimensional Euclidean space at a given time $t>0$. Let $\Omega_0$ be the reference configuration, which we assume to coincide with the initial configuration, and be $\mathbf{X}\in\Omega_0$ a generic point. The motion of the body $\mathcal{B}$ is given by \vspace{-0.1cm}
\begin{equation*}
\chi(\mathbf{X},t;\bm\mu) = \mathbf{x} \quad \forall t>0, \vspace{-0.1cm}
\end{equation*}
which maps the material position $\mathbf{X}$ in the reference configuration $\Omega_0$ to the spatial position $\mathbf{x}$ in the deformed or current configuration $\Omega_t$ for all times $t>0$. A motion $\chi$ of a body $\mathcal{B}$ will change the body's shape, position and/or orientation. For a given parameter vector $\bm\mu\in\mathcal{P}$, the displacement vector field \vspace{-0.1cm}
\begin{equation*}
\mathbf{u}(\mathbf{X},t;\bm\mu) = \chi(\mathbf{X},t;\bm\mu) - \mathbf{X} \vspace{-0.1cm}
\end{equation*}
relates the position $\mathbf{X}$ of a particle in the reference configuration to its position $\mathbf{x}$ in the current configuration at time $t>0$. A crucial quantity in nonlinear mechanics is the deformation gradient \vspace{-0.1cm}
\begin{equation*}
\mathbf{F}(\mathbf{X},t;\bm\mu) = \frac{\partial \chi(\mathbf{X},t;\bm\mu)}{\partial \mathbf{X}} = \mathbf{I} + \nabla_0 \mathbf{u}(\mathbf{X},t;\bm\mu), \vspace{-0.1cm}
\end{equation*}
which characterizes changes of material elements during motion. The change in volume between the reference and the current configurations at time $t>0$ is given by $J(\mathbf{X},t;\bm\mu) = \det\mathbf{F}(\mathbf{X},t;\bm\mu)>0$. Common measures of strain are the right Cauchy-Green strain tensor and the Green-Lagrange strain tensor, that are defined as \vspace{-0.1cm}
\begin{equation}\label{eq:strain}
\mathbf{C} = \mathbf{F}^T\mathbf{F}, \qquad \mathbf{E} = \frac{1}{2}(\mathbf{C}-\mathbf{I}), \vspace{-0.1cm}
\end{equation}
respectively. The equation of motion for a continuous medium is given by the conservation of mass and the balance of the linear momentum, 
in material coordinates, reads as follows: 
\begin{equation*}
\rho_0\partial_t^2\mathbf{u}(\mathbf{X},t;\bm\mu) - \nabla_0\cdot\mathbf{P}(\mathbf{F}(\mathbf{X},t;\bm\mu)) = \mathbf{b}_0(\mathbf{X},t;\bm\mu), \qquad \mathbf{X}\in\Omega_0,~t>0
\end{equation*}
where $\rho_0$ is the density of the body, $\mathbf{P(F)}$ is the first Piola-Kirchhoff stress tensor and $\mathbf{b}_0$ is an external body force. Proper boundary and initial conditions must be specified to ensure the well-posedness of the problem. In addition, we need a constitutive equation for $\mathbf{P}$, that is, a stress-strain relationship describing the material behavior. Here, we consider hyperelastic materials, for which  the existence of a strain density function $\mathcal{W}\colon Lin^+\rightarrow\mathbb{R}$ such that
\begin{equation*}
\mathbf{P(F)} = \frac{\partial\mathcal{W}(\mathbf{F})}{\partial\mathbf{F}}
\end{equation*}
 is postulated. Note that, since $\mathbf{F}$ depends on the displacement $\mathbf{u}$, we can equivalently write $\mathbf{P(F)}$ or $\mathbf{P}(\mathbf{u})$. The strong formulation of a general initial boundary-valued problem in elastodynamics thus reads as follows: given a body force $\mathbf{b}_0 = \mathbf{b}_0(\mathbf{X},t;\bm\mu)$, a prescribed displacement $\bar{\mathbf{u}} = \bar{\mathbf{u}}(\mathbf{X},t;\bm\mu)$ and surface traction $\bar{\mathbf{T}} = \bar{\mathbf{T}}(\mathbf{X},t,\mathbf{N};\bm\mu)$, find the unknown displacement field $\mathbf{u}(\bm\mu)\colon\Omega_0\times(0,T]\rightarrow\mathbb{R}^3$ so that
\begin{align}\label{eq:strong-IBVP}
\left\{ \begin{array}{lllr}
\rho_0\partial^2_t {\mathbf{u}}(\mathbf{X},t;\bm\mu) - \nabla_0\cdot\mathbf{P}(\mathbf{u}(\mathbf{X},t;\bm\mu)) = \mathbf{b}_0(\mathbf{X},t;\bm\mu) && \text{in } & \Omega_0\times(0,T]\\
\mathbf{u}(\mathbf{X},t;\bm\mu) = \bar{\mathbf{u}}(\mathbf{X},t;\bm\mu) && \text{on } & \Gamma_0^{D}\times(0,T]\\
\mathbf{P}(\mathbf{u}(\mathbf{X},t;\bm\mu))\mathbf{N} = \bar{\mathbf{T}}(\mathbf{X},t,\mathbf{N};\bm\mu) && \text{on } & \Gamma_0^{N}\times(0,T]\\
\mathbf{P}(\mathbf{u}(\mathbf{X},t;\bm\mu))\mathbf{N} + \alpha \mathbf{u}(\mathbf{X},t;\bm\mu) + \beta \partial_t\mathbf{u}(\mathbf{X},t;\bm\mu) = \mathbf{0} && \text{on } & \Gamma_0^{R}\times(0,T]\\
\mathbf{u}(\mathbf{X},0;\bm\mu) = \mathbf{u}_0(\mathbf{X};\bm\mu),~~ \partial_t{\mathbf{u}}(\mathbf{X},0;\bm\mu) = \dot{\mathbf{u}}_0(\mathbf{X};\bm\mu) && \text{in } & \Omega_0\times\{0\}
\end{array} \right.
\end{align}
where $\mathbf{N}$ is the outer normal unit vector and $\alpha,\beta\in\mathbb{R}$. The boundary of the reference domain is divided such that $\Gamma_0^D\cup\Gamma_0^N\cup\Gamma_0^R = \Gamma$, with $\Gamma_0^i\cap\Gamma_0^j=\emptyset$ for $i,j\in\{D,N,R\}$.  The corresponding variational form can we written as: $\forall t\in(0,T]$, find the unknown displacement field $\mathbf{u}(t;\bm\mu)\in V$ such that
\begin{align}\label{eq:weak-IBVP}
\langle R(\mathbf{u}(t;\bm\mu),t;\bm\mu),\bm\eta\rangle &:= \int_{\Omega_0} \rho_0\partial_t^2{\mathbf{u}}(t;\bm\mu) \cdot{\bm\eta}d\Omega + \int_{\Omega_0} \mathbf{P}(\mathbf{u}(t;\bm\mu)):\nabla{\bm\eta}d\Omega \nonumber  \\
&+ \int_{\Gamma_0^R} \left(\alpha\mathbf{u}(t;\bm\mu)+\beta\partial_t\mathbf{u}(t;\bm\mu)\right)\cdot{\bm\eta}d\Gamma - \int_{\Gamma_N} \bar{\mathbf{T}}(t,\mathbf{N};\bm\mu)\cdot{\bm\eta}d\Gamma - \int_{\Omega_0} \mathbf{b}_0(t;\bm\mu)\cdot{\bm\eta}d\Omega = 0\\
\int_{\Omega_0} \mathbf{u}(0;\bm\mu) \cdot{\bm\eta}d\Omega &= \int_{\Omega_0} \mathbf{u}_0(\bm\mu) \cdot{\bm\eta}d\Omega \nonumber, \qquad \int_{\Omega_0} \partial_t\mathbf{u}(0;\bm\mu) \cdot{\bm\eta}d\Omega = \int_{\Omega_0} \dot{\mathbf{u}}_0(\bm\mu) \cdot{\bm\eta}d\Omega \nonumber
\end{align}
for any test function $\bm\eta$, where $V=V(\Omega_0)$ denotes a suitable Hilbert space on the reference configuration $\Omega_0\in\mathbb{R}^3$ and $V'$ its dual. This equation is inherently nonlinear and additional source of nonlinearity is introduced in the material law, i.e. when using a nonlinear $\mathcal{W} =\mathcal{W}(\mathbf{F})$, which is often the case of engineering applications. 

For the sake of simplicity, in all test cases, we neglect the body forces $\mathbf{b}_0(\bm\mu)$ and consider zero initial conditions $\mathbf{u}_0(\bm\mu) = \dot{\mathbf{u}}_0(\bm\mu) = \mathbf{0}$. Regarding boundary conditions, we consider $\bar{\mathbf{u}}(\bm\mu) = \mathbf{0}$ on the Dirichlet boundary $\Gamma_0^D$ and always assume $\alpha=\beta=0$, so that we actually impose homogeneous Neumann conditions on $\Gamma_0^R$. Finally, the traction vector is given by 
\begin{equation*}
\bar{\mathbf{T}}(\mathbf{X},t,\mathbf{N};\bm\mu) = -\mathbf{g}(t;\bm\mu)J\mathbf{F}^{-T}\mathbf{N},
\end{equation*}
where $\mathbf{g}(t;\bm\mu)$ represents an external load and will be   specified according to the application at hand.

The residual in (\ref{eq:residual}) is given by
\begin{equation*}
\begin{aligned}
\mathbf{R}(\mathbf{u}_h^{n}(\bm\mu), t^{n};\bm\mu) &:= \left(\dfrac{\rho_0}{\Delta t^2}\mathcal{M} + \dfrac{1}{\Delta t}\mathcal{F}_{\beta}^{int} + \mathcal{F}_{\alpha}^{int}\right)\mathbf{u}_h^{n}(\bm\mu) + \mathcal{S}(\mathbf{u}_h^{n}(\bm\mu)) - \left(\dfrac{2\rho_0}{\Delta t^2}\mathcal{M} + \dfrac{1}{\Delta t}\mathcal{F}_{\beta}^{int}\right)\mathbf{u}_h^{n-1}(\bm\mu) \\ & + \dfrac{\rho_0}{\Delta t^2}\mathcal{M}\mathbf{u}_h^{n-2}(\bm\mu) - \mathcal{F}^{ext,n}(\bm\mu),
\end{aligned}
\end{equation*}
for $n=1,\dots,N_t$, where $\mathbf{u}_h^0(\bm\mu)$ and $\mathbf{u}_h^{-1}(\bm\mu)$ are known for the initial condition, and
\begin{align*}
&[\mathcal{M}]_{ij} = \int_{\Omega_0} \bm\varphi_j\cdot\bm\varphi_id\Omega, 
\qquad 
[\mathcal{F}_{\beta}^{int}]_{ij} = \int_{\Gamma_0^R} \beta~\bm\varphi_j\cdot\bm\varphi_i d\Gamma, \qquad
[\mathcal{F}_{\alpha}^{int}]_{ij} = \int_{\Gamma_0^R} \alpha~\bm\varphi_j\cdot\bm\varphi_i d\Gamma,\\
&[\mathcal{S}(\mathbf{u}_h^n(\bm\mu))]_{i} = \int_{\Omega_0} \mathbf{P}(\mathbf{u}_h^n(\bm\mu))\colon \nabla{\bm\varphi_i}d\Omega, 
\qquad 
[\mathcal{F}^{ext,n}(\bm\mu)]_{i} = \int_{\Gamma_0^N} \bar{\mathbf{T}}^n(\mathbf{N};\bm\mu)\cdot{\bm\varphi_i}d\Gamma + \int_{\Omega_0} \mathbf{b}_0^n(\bm\mu)\cdot{\bm\varphi_i}d\Omega,
\end{align*}
for all $i,j=1,\dots,N_h$, being $\{\bm\varphi_i\}_{i=1}^{N_h}$ a basis for the finite element (FE) space. 

\newpage

As a measure of accuracy of the reduced approximations with respect to the FOM solution, we consider time-averaged $L^2$-errors of the displacement vector, that are defined as follows: \vspace{-0.2cm}
\begin{equation}\label{eq:error}
\begin{aligned}
\epsilon_{abs}(\bm\mu) &= \frac{1}{N_t}\sum_{n=1}^{N_t} \lVert \mathbf{u}_h(\cdot,t^n;\bm\mu) - \mathbf{Vu}_N(\cdot,t^n;\bm\mu)\rVert_{2}, \vspace{-0.1cm}\\
\epsilon_{rel}(\bm\mu) &= \frac{1}{N_t}\sum_{n=1}^{N_t} \frac{\lVert \mathbf{u}_h(\cdot,t^n;\bm\mu) - \mathbf{Vu}_N(\cdot,t^n;\bm\mu)\rVert_{2}}{\lVert \mathbf{u}_h(\cdot,t^n;\bm\mu)\rVert_{2}}.
\end{aligned}
\end{equation}
The CPU time ratio, that is the ratio between FOM and ROM computational times, is used to measure efficiency, since it represents the speed-up offered by  the ROM with respect to the FOM. The code is implemented in Python in our software package pyfe$^\text{x}$, a Python binding with the in-house Finite Element library \texttt{life$^\texttt{x}$} (\url{https://lifex.gitlab.io/lifex}), a high-performance C++ library based on the \texttt{deal.II} (\url{https://www.dealii.org}) Finite Element core \cite{dealII92}. Computations have been performed on a PC desktop computer with 3.70GHz Intel Core i5-9600K CPU and 16GB RAM.


\subsection{Deformation of a clamped rectangular beam}\label{sec:DB}

\begin{figure}[b!]
	\centering
	\includegraphics[width=0.795\textwidth]{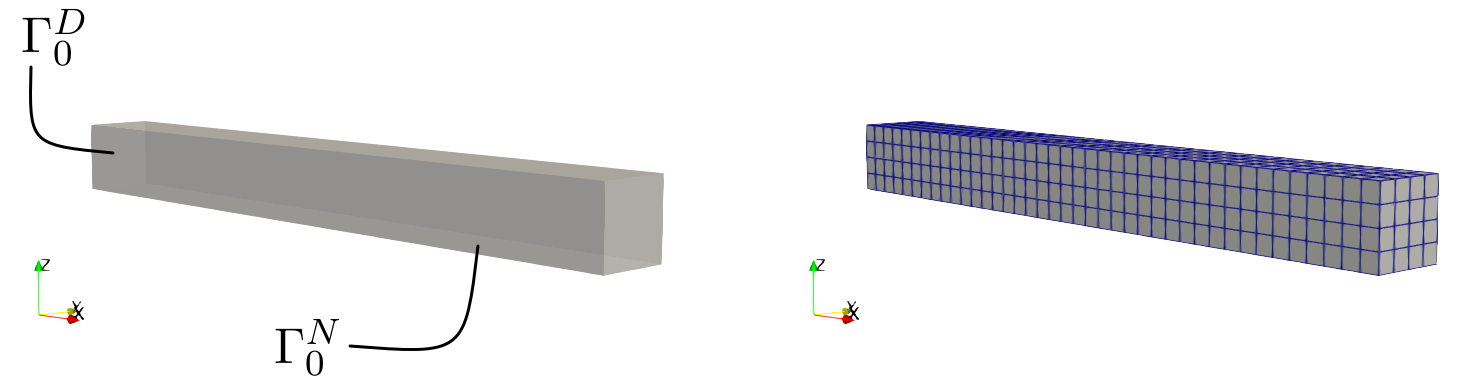}
	\caption{Rectangular beam geometry (left) and computational grid (right).}
	\label{fig:DB_mesh}
\end{figure}

The first series of test cases represents a typical structural mechanical problem, with reference geometry $\bar\Omega_0 = [0,10^{-2}]\times[0,10^{-3}]\times[0,10^{-3}]$ m$^3$, reported in Figure~\ref{fig:DB_mesh}.
For the continuum body $\mathcal{B}$ under investigation, we consider a nearly-incompressible neo-Hookean material, which is characterized by the following strain density energy function \vspace{-0.1cm}
\begin{equation*}
\mathcal{W}(\mathbf{F}) = \frac{G}{2}(\mathcal{I}_1 - 3) + \frac{K}{4}( (J-1)^2 + \ln^2(J) ), \vspace{-0.1cm}
\end{equation*}
where $G>0$ is the shear modulus, $\mathcal{I}_1 = J^{-\frac{2}{3}}\det(\mathbf{C})$ and the latter term is needed to enforce in\-com\-pres\-si\-bi\-li\-ty, being the bulk modulus $K>0$ the penalization factor. This choice leads to the following first Piola-Kirchhoff stress tensor, characterized by a nonpolynomial nonlinearity, \vspace{-0.1cm}
\begin{equation*}
\mathbf{P}(\mathbf{F}) = GJ^{-\frac{2}{3}}\left(\mathbf{F} - \frac{1}{3}\mathcal{I}_1\mathbf{F}^T\right) + \frac{K}{2}J\left(J-1+\frac{1}{J}\ln(J)\right)\mathbf{F}^T. \vspace{-0.1cm}
\end{equation*}
The beam is clamped at the left-hand side, that is, Dirichlet boundary conditions are imposed on the left face $x=0$, whilst a pressure load changing with the deformed surface orientation is applied to the entire bottom face $z=0$ (i.e. $\Gamma_0^N$). Homogeneous Neumann conditions are applied on the remaining boundaries (i.e. $\Gamma_0^R$ with $\alpha=\beta=0$). As possible functions for the external load $\mathbf{g}(t;\bm\mu)$, we choose
\begin{enumerate}
	\item a linear function  $\mathbf{g}(t;\bm\mu) = \widetilde{p}~t/T$;
	\item a triangular or hat function $\mathbf{g}(t;\bm\mu) = \widetilde{p}~\left(2t~\chi(t)_{\left(0,\frac{T}{2}\right]} + 2(T-t)~\chi(t)_{\left(\frac{T}{2},T\right]}\right)$;
	\item a step function $\mathbf{g}(t;\bm\mu) = \widetilde{p}~\chi(t)_{\left(0,\frac{T}{3}\right]}$, so that the presence of the inertial term is not negligible.	
\end{enumerate}
Here, $\widetilde{p}>0$ is a parameter controlling the maximum load. The FOM is built on a hexahedral mesh with $640$ elements and $1025$ vertices, resulting in a high-fidelity dimension $N_h=3075$ (since $\mathbb{Q}_1$-FE are employed). The resulting computational mesh in the reference configuration is reported in Figure~\ref{fig:DB_mesh}.
\newpage

The following sections are organized as follows: first, we analyze the accuracy and the efficiency of the ROM without hyper-reduction with respect to the POD tolerance $\varepsilon_{POD}$, thus resulting in reduced subspaces of different dimensions $N\in\mathbb{N}$. Then, for a fixed basis $\mathbf{V}\in\mathbb{R}^{N_h\times N}$, POD-Galerkin-DEIM approximation capabilities are investigated for different sizes of the reduced mesh, associated with different values of the tolerance $\varepsilon_{DEIM}$ for the computation of the residual basis $\bm\Phi_{\mathcal{R}}\in\mathbb{R}^{N_h\times m}$. Finally, the performances of Deep-HyROMnet are assessed and compared to those of DEIM-based hyper-ROMs.


\subsubsection{Test case 1: linear function for the pressure load}

Let us consider the parametrized linear function 
\begin{equation*}
\mathbf{g}(t;\bm\mu) = \widetilde{p}~t/T, 
\end{equation*}
for the pressure load, describing a situation in which a structure is progressively loaded. We choose a time interval $t\in[0,0.25]$~s and employ a uniform time step $\Delta t = 5\cdot10^{-3}$~s for the time discretization scheme, resulting in a total number of $50$ time iterations. As parameters, we consider:
\begin{itemize}
	\item the shear modulus $G\in[0.5\cdot10^4,1.5\cdot10^4]$~Pa;
	\item the bulk modulus $K\in[2.5\cdot10^4,7.5\cdot10^4]$~Pa;
	\item the external load parameter $\widetilde{p}\in[2,6]$~Pa.
\end{itemize}

\begin{figure}[b!]
	\centering
	\includegraphics[width=0.45\textwidth]{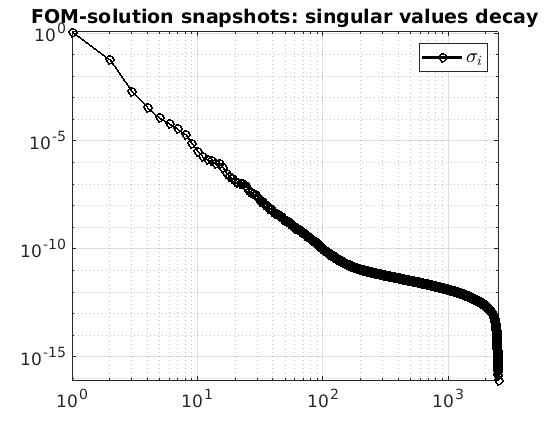}
	\caption{Test case 1. Decay of the singular values of the FOM solution snapshots matrix.}
	\label{fig:DB_ramp_svd_uh}
\end{figure}

Given a training set of $n_s=50$ points generated from the three-dimensional parameter space $\mathcal{P}$ through latin hypercube sampling (LHS), we compute the reduced basis $\mathbf{V}\in\mathbb{R}^{N_h\times N}$ using the POD method with tolerance \vspace{-0.2cm}
\begin{equation*}
\varepsilon_{POD}\in\{10^{-3}, 5\cdot10^{-4}, 10^{-4}, 5\cdot10^{-5}, 10^{-5}, 5\cdot10^{-6}, 10^{-6}\}. \vspace{-0.2cm}
\end{equation*}
The corresponding reduced dimensions are $N=3$, $4$, $5$, $6$, $8$, $9$ and $15$, respectively. In Figure~\ref{fig:DB_ramp_svd_uh} we show the singular values of the snapshot matrix related to the FOM displacement $\mathbf{u}_h$, where a rapid decay of the plotted quantity means that a small number of RB functions are needed to correctly approximate the high-fidelity solution.

\begin{figure}[t!]
\vspace{-0.2cm}
	\centering
	\begin{subfigure}{0.49\textwidth}
		\centering
		\includegraphics[width=0.95\textwidth]{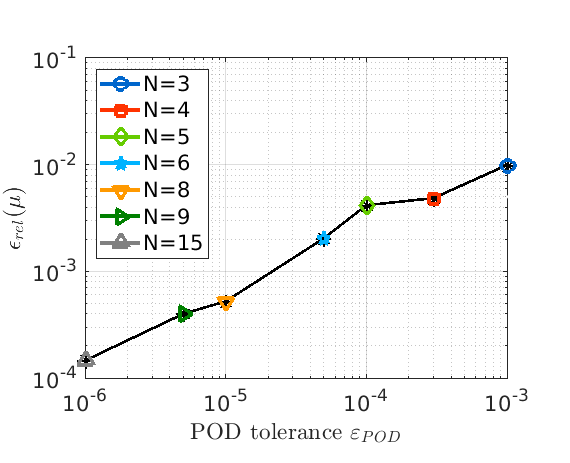}
	\end{subfigure}
	\begin{subfigure}{0.49\textwidth}
		\centering
		\includegraphics[width=0.95\textwidth]{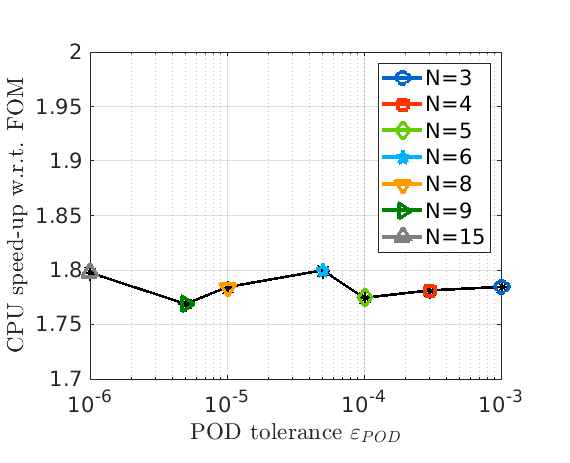}
	\end{subfigure}
	\caption{Test case 1. Average over 50 testing parameters of relative error $\epsilon_{rel}$ (left) and average speed-up (right) of ROM without hyper-reduction.}
	\label{fig:DB_ramp_ROM}
\end{figure}

The average relative error $\epsilon_{rel}$ between the FOM and the POD-Galerkin ROM solutions computed over a testing set of 50 randomly chosen parameters, different from the ones used to compute the solution snapshots, is reported in Figure~\ref{fig:DB_ramp_ROM}, together with the CPU time ratio. The approximation error decreases up to an order of magnitude when reducing the POD tolerance $\varepsilon_{POD}$ from $10^{-3}$ to $10^{-6}$, corresponding to an increase of the RB dimension from $N=3$ to $N=15$. Despite being the RB space low-dimensional, the computational speed-up achieved by the reduced model is negligible. This is due to the fact that the ROM still depends on the FOM dimension $N_h$ during the online stage. For this reason, we need to rely on suitable hyper-reduction techniques.

For the construction of both hyper-reduced models (POD-Galerkin-DEIM and Deep-HyROMnet), we need first to compute snapshots from the ROM solutions for given parameter values and time instants, in order to build either the DEIM basis $\bm\Phi_\mathcal{R}$ or train the DNNs $\bm\rho_N$ and $\bm\iota_N$. To this goal, we choose a POD-Galerkin ROM with dimension $N=4$, being it a good balance between accuracy and computational effort for the test case at hand, and perform ROM simulations for a given set of $n_s'=200$ parameter samples to collect residual and Jacobian data.

In order to investigate the impact of hyper-reduction onto the ROM solution reconstruction error, we compute the DEIM basis $\mathbf{\bm\Phi}_{\mathcal{R}}$ for the approximation of the residual using the POD method with different tolerances, that are
\begin{equation*}
\varepsilon_{DEIM}\in\{10^{-3}, 5\cdot10^{-4}, 10^{-4}, 5\cdot10^{-5}, 10^{-5}, 5\cdot10^{-6}, 10^{-6}\},
\end{equation*}
corresponding to $m=22$, $25$, $30$, $33$, $39$, $43$, $51$, respectively. Larger POD tolerances were not sufficient to ensure the convergence of Newton method for all considered combinations of parameters, so that higher speed-ups cannot be achieved by decreasing the basis dimension $m$. 

\begin{figure}[b!]
\vspace{-0.2cm}
	\centering
	\begin{subfigure}{0.49\textwidth}
		\centering
		\includegraphics[width=0.95\textwidth]{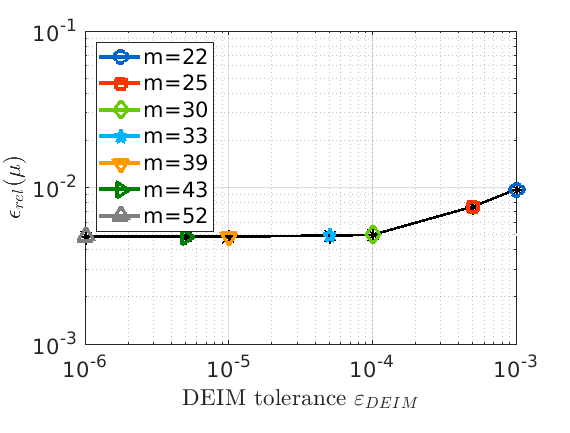}
	\end{subfigure}
	\begin{subfigure}{0.49\textwidth}
		\centering
		\includegraphics[width=0.95\textwidth]{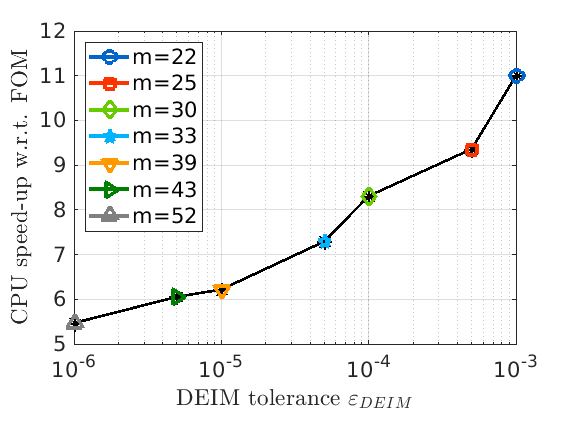}
	\end{subfigure}
	\caption{Test case 1. Average over 50 testing parameters of relative error $\epsilon_{rel}$ (left) and average speed-up (right) of POD-Galerkin-DEIM with $N=4$.}
	\label{fig:DB_ramp_DEIM_N4}
\end{figure}
The average relative error $\epsilon_{rel}$ is evaluated over the testing set and plotted in Figure~\ref{fig:DB_ramp_DEIM_N4}, as well as the CPU time ratio. To compute the high-fidelity solutions, $26$~s are required in average, while a POD-Galerkin-DEIM ROM, with $N=4$ and $m=22$, requires only $2.4$~s, thus yielding a speed-up of $\times11$ compared to the FOM.

Data related to the performances of the POD-Galerkin-DEIM method for $N=4$ and different values of $m$ are shown in Table~\ref{tab:DB_ramp}. The number of elements of the reduced mesh represents a small percentage of the one forming the original grid, so that the cost related to the residual assembling is remarkably alleviated. Nonetheless, it is obvious that the main computational bottleneck is the construction of the reduced system at each Newton iteration, and in particular the assembling of the residual vector on the reduced mesh, which requires between $78\%$ and $88\%$ of the total (online) CPU time. In particular, almost $90\%$ of this computational time is demanded for assembling the residual $\mathbf{R}(\mathbf{Vu}_N^{n,(k)}(\bm\mu),t^{n};\bm\mu)$ on the reduced mesh, while computing the associated Jacobian matrix using the automatic differentiation tool takes less than $1\%$.
\begin{table}[h!]
	\centering
	\begin{tabular}{|l||c|c|c|} 
		\hline
		POD tolerance $\varepsilon_{DEIM}$ & $5\cdot10^{-4}$ & $5\cdot10^{-5}$ & $5\cdot10^{-6}$\\
		\hline
		DEIM interpolation dofs $m$ & $25$ & $33$ & $43$ \\
		\hline
		Reduced mesh elements (total: $640$) & $86$ & $115$ & $168$ \\
		\hline
		\hline
		Online CPU time & $2.8$~s & $3.6$~s & $4.3$~s \\
		$\quad\circ$ system construction $[*]$ & $78\%$ & $83\%$ & $88\%$\\
		$\quad\circ$ system solution   & $0.16\%$ & $0.13\%$ & $0.09\%$\\
		\hline
		\hline
		$[*]$ System construction for each Newton iteration & $0.02$~s & $0.02$~s & $0.03$~s \\
		$\quad\circ$ residual assembling              & $89\%$  & $87\%$  & $88\%$\\
		$\quad\circ$ Jacobian computing through AD    & $0.6\%$ & $0.4\%$ & $0.5\%$\\
		\hline
		\hline
		Computational speed-up & $\times$9.4 & $\times$7.3 & $\times$6.0 \\
		\hline
		Time-averaged $L^2(\Omega_0)$-absolute error & $3\cdot10^{-5}$ & $2\cdot10^{-5}$ & $2\cdot10^{-5}$\\
		\hline
		Time-averaged $L^2(\Omega_0)$-relative error & $8\cdot10^{-3}$ & $5\cdot10^{-3}$ & $5\cdot10^{-3}$ \\
		\hline
	\end{tabular}
	\caption{Test case 1. Computational data related to POD-Galerkin-DEIM with $N=4$ and different values of $m$.}
	\label{tab:DB_ramp}
\end{table}

Finally, we analyze the performances of Deep-HyROMnet and compare them in terms of both accuracy and efficiency with POD-Galerkin-DEIM ROMs. The average of the absolute error $\epsilon_{abs}$, the relative error $\epsilon_{rel}$ and the CPU time ratio are reported in Table~\ref{tab:DB_ramp_N4_hyper-ROM}. In terms of efficiency, the DNN-based ROMs outperform the DEIM-based hyper-ROMs substantially, being almost $100$ times faster than POD-Galerkin-DEIM ROM with $m=22$, whist achieving the same accuracy. In particular, Deep-HyROMnet is able to compute the reduced solutions in less than 0.03~s, thus yielding an overall speed-up of order $\mathcal{O}(10^3)$ compared to the FOM.
\begin{table}[b!]
	\centering
	\vspace{0.15cm}
	\begin{tabular}{|l||c|c|c|} 
		\hline
		& DEIM ($m=$22) & DEIM ($m=$30) & Deep-HyROMnet\\
		\hline
		\hline
		Computational speed-up & $\times$11  & $\times$8 & $\times$1012\\
		\hline
		Avg. CPU time & 2~s & 3~s & 0.026~s\\
		\hline
		Time-avg. $L^2(\Omega_0)$-absolute error & $7.4\cdot10^{-5}$ & $1.9\cdot10^{-5}$ & $7.7\cdot10^{-5}$\\
		\hline
		Time-avg. $L^2(\Omega_0)$-relative error & $9.7\cdot10^{-3}$ & $5.0\cdot10^{-3}$ & $8.3\cdot10^{-3}$\\
		\hline
	\end{tabular}
	\caption{Test case 1. Computational data related to POD-Galerkin-DEIM ROMs and Deep-HyROMnet, for $N=4$.}
	\label{tab:DB_ramp_N4_hyper-ROM} 
\end{table}

\begin{figure}[t!]
	\centering
	\includegraphics[width=0.95\textwidth]{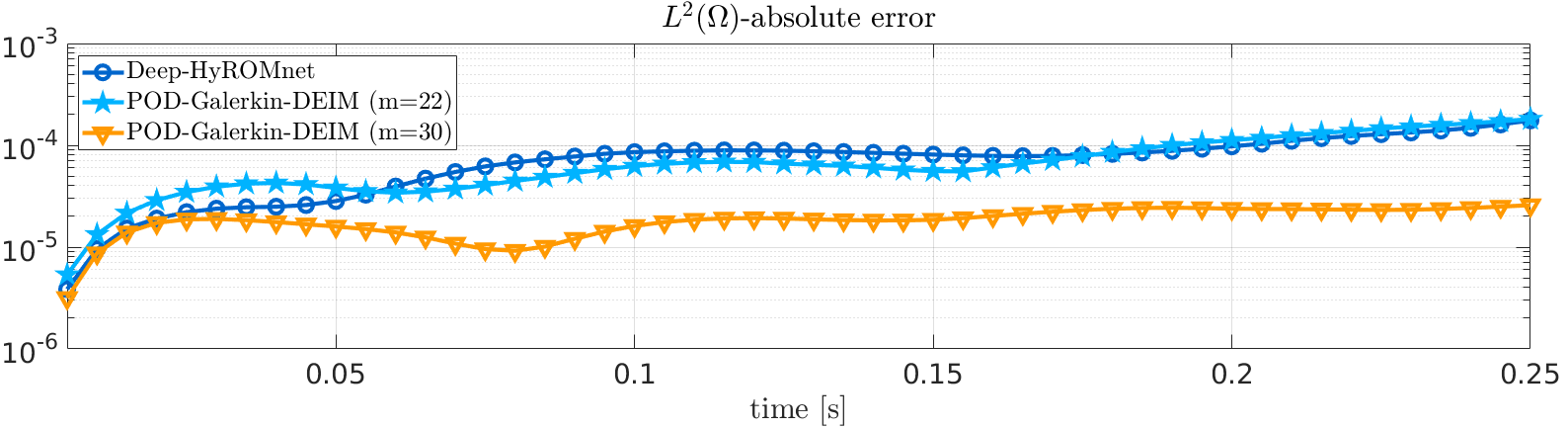}
	\caption{Test case 1. Evolution in time of the average $L^2(\Omega_0)$-absolute error for $N=4$ computed using POD-Galerkin-DEIM and Deep-HyROMnet.}
	\label{fig:DB_ramp_N4_hyper-ROM_err_abs}
\end{figure}

\begin{figure}[t!]
	\centering
	\includegraphics[width=0.8\textwidth]{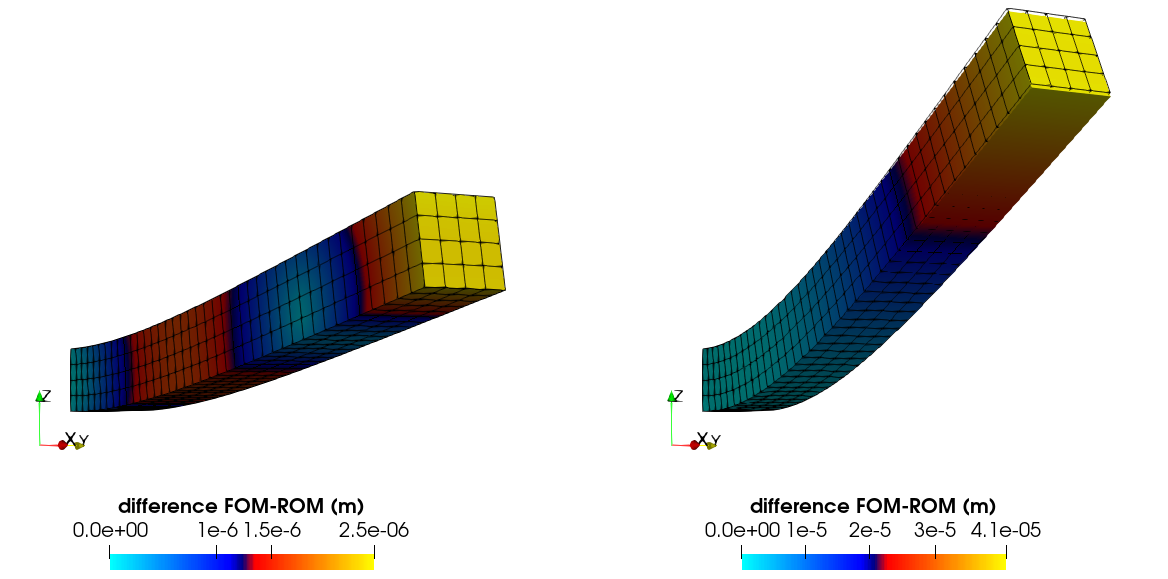}
	\caption{Test case 1. FOM (wireframe) and Deep-HyROMnet (colored) solutions at time $T=0.25$~s for $\bm\mu = [1.3225\cdot10^4~\text{Pa}, 3.9875\cdot10^4~\text{Pa}, 3.43~\text{Pa}]$ (left) and $\bm\mu = [0.6625\cdot10^4~\text{Pa}, 5.8625\cdot10^4~\text{Pa}, 4.89~\text{Pa}]$ (right).}
	\label{fig:DB_ramp_N4_DeepHyROMnet_mu10_mu13_error}
\end{figure}

The evolution of the $L^2(\Omega_0)$-absolute error, averaged over the testing parameters, is reported in Figure~\ref{fig:DB_ramp_N4_hyper-ROM_err_abs} for all of the hyper-ROMs considered.
The final accuracy of the hyper-ROMs equals that of the ROM without hyper-reduction, i.e. $\epsilon_{rel}\approx10^{-2}$, meaning that the projection error dominates over the nonlinear operators approximation error. The difference between the FOM and Deep-HyROMnet solutions at time $T=0.25$~s is shown in Figure \ref{fig:DB_ramp_N4_DeepHyROMnet_mu10_mu13_error} in two scenarios.

In order to increase the accuracy of the reduced solution, we should consider higher values of the RB dimension $N$. As a matter of fact, by increasing the RB dimension, the task of the DNNs becomes more complex, meaning that more training samples and a larger size of the neural networks themselves may be required. In Table~\ref{tab:DB_ramp_N8_hyper-ROM} are reported the computational data associated with POD-Galerkin-DEIM and Deep-HyROMnet hyper-ROMs when $N=8$, but the same number of training snapshots and the same DNN architectures of previous case (i.e. $N=4$) are employed. We observe that POD-Galerkin-DEIM is able to provide more accurate approximations of the high-fidelity solution by increasing the size $m$ of the residual basis, albeit reducing the online speed-up with respect to the FOM. On the other hand, in the context of multi-query problems, such as uncertainty quantification or optimization, where thousands of queries to the parameter-to-solution map are required, it is of paramount importance to decrease the CPU time needed for the solution of the reduced problem at each new instance of the input parameter vector.
\begin{table}[h]
	\centering
	\begin{tabular}{|l||c|c|c|} 
		\hline
		& DEIM ($m=$29) & DEIM ($m=$51) & Deep-HyROMnet\\
		\hline
		\hline
		Computational speed-up & $\times$8  & $\times$5 & $\times$949\\
		\hline
		Avg. CPU time & 3~s & 5~s & 0.027~s\\
		\hline
		Time-avg. $L^2(\Omega_0)$-absolute error & $2.6\cdot10^{-5}$ & $3.9\cdot10^{-6}$ & $9.0\cdot10^{-5}$\\
		\hline
		Time-avg. $L^2(\Omega_0)$-relative error & $1.1\cdot10^{-2}$ & $6.0\cdot10^{-4}$ & $8.1\cdot10^{-3}$\\
		\hline
	\end{tabular}
	\caption{Test case 1. Computational data related to POD-Galerkin-DEIM ROMs and Deep-HyROMnet, for $N=8$.}
	\label{tab:DB_ramp_N8_hyper-ROM}
	\vspace{-0.2cm}
\end{table}
 

\subsubsection{Test case 2: hat function for the pressure load}

\begin{figure}[h]
	\centering
	\includegraphics[width=0.45\textwidth]{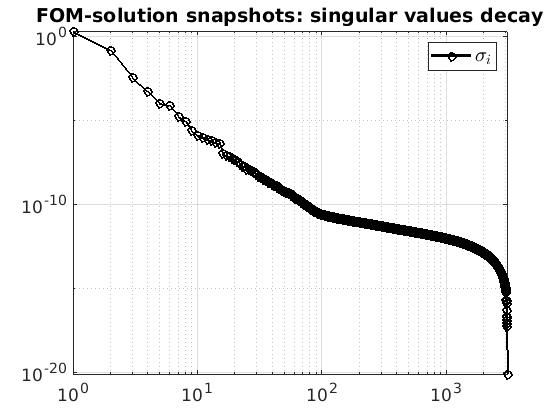}
	\caption{Test case 2. Decay of the singular values of the FOM solution.}
	\label{fig:DB_hat_svd_uh}
\end{figure}

Let us now consider a piecewise linear pressure load given by the following hat function 
\begin{equation*}
\mathbf{g}(t;\bm\mu) = \widetilde{p}~\left(2t~\chi(t)_{\left(0,\frac{T}{2}\right]} + 2(T-t)~\chi(t)_{\left(\frac{T}{2},T\right]}\right),
\end{equation*}
describing the case in which a structure is increasingly loaded until a maximum pressure is reached, and then linearly unloaded in order to recover the initial resting state. For the case at hand, we choose $t\in[0,0.35]$~s and $\Delta t = 5\cdot10^{-3}$~s, resulting in a total number of $70$ time steps. As parameter, we consider the external load parameter $\widetilde{p}\in[2,12]$ Pa. The shear modulus $G$ and the bulk modulus $K$ are fixed to the values $10^4$ Pa and $5\cdot10^4$ Pa, respectively. Let us consider a training set of $n_s=50$ points generated from the one-dimensional parameter space $\mathcal{P}=[2,12]$ Pa through LHS and build the RB basis $\mathbf{V}\in\mathbb{R}^{N_h\times N}$ with $N=4$, corresponding to $\varepsilon_{POD}=10^{-4}$. The singular values of the solution snapshots matrix are reported in Figure~\ref{fig:DB_hat_svd_uh}.  

Given the Galerkin-ROM nonlinear data collected for $n_s'=300$ sampled parameters, the DEIM residual basis $\mathbf{\bm\Phi}_{\mathcal{R}}$ is computed using the POD method with tolerance  \vspace{-0.1cm}
\begin{equation*}
\varepsilon_{DEIM}\in\{10^{-3}, 5\cdot10^{-4}, 10^{-4}, 5\cdot10^{-5}, 10^{-5}, 5\cdot10^{-6}, 10^{-6}\}, \vspace{-0.1cm}
\end{equation*}
corresponding to $m=14$, $16$, $22$, $24$, $29$, $31$, $37$, respectively. Tolerances $\varepsilon_{DEIM}$ larger than the values reported above were not sufficient to ensure convergence of Newton method for all the considered parameters. 
\begin{figure}[h]
	\centering
	\begin{subfigure}{0.49\textwidth}
		\centering
		\includegraphics[width=0.95\textwidth]{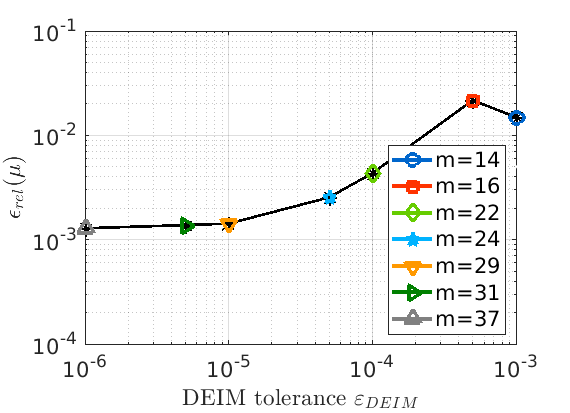}
	\end{subfigure}
	\begin{subfigure}{0.49\textwidth}
		\centering
		\includegraphics[width=0.95\textwidth]{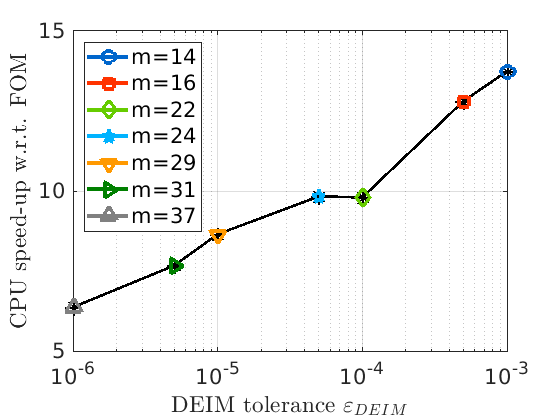}
	\end{subfigure}
	\caption{Test case 2. Average over 50 testing parameters of relative error $\epsilon_{rel}$ (left) and average speed-up (right) of POD-Galerkin-DEIM with $N=4$.}
	\label{fig:DB_hat_DEIM_N4}
\end{figure}
The relative error $\epsilon_{rel}$, evaluated over a testing set of $50$ parameters, is about $10^{-2}$ when using $m=14$ residual basis, and can be further reduced of one order of magnitude when increasing the DEIM dimension to $m=29$, albeit highly decreasing the CPU time ratio, as shown in Figure~\ref{fig:DB_hat_DEIM_N4}.

Table \ref{tab:DB_hat_hyper-ROM_N4} shows the comparison between POD-Galerkin-DEIM (with $m=14$ and $m=29$) and Deep-HyROMnet hyper-reduced models on a testing set of 50 parameter instances. As observed in the previous test case, Deep-HyROMnet is able to achieve good results in terms of accuracy, comparable with the fastest DEIM-based model ($m=14$), at a greatly reduced cost. Also in this case, the speed-up achieved by our DNN-based hyper-ROM is of order $\mathcal{O}(10^3)$ with respect to the FOM, since less than $0.04$~s are needed to compute the reduced solution for each new instance of the parameter, against a time of about $40$~s required by the FOM, and of $3$~s required by POD-Galerkin-DEIM. 
\begin{table}[t]
	\centering
	\begin{tabular}{|l||c|c|c|} 
		\hline
		& DEIM ($m=14$) & DEIM ($m=29$) & Deep-HyROMnet \\
		\hline
		\hline
		Computational speed-up & $\times$14  & $\times$9 & $\times$1153\\
		\hline
		Avg. CPU time & 3~s & 5~s & 0.035~s\\
		\hline
		Time-avg. $L^2(\Omega_0)$-absolute error & $9.0\cdot10^{-5}$ & $6.8\cdot10^{-6}$ & $2.0\cdot10^{-4}$\\
		\hline
		Time-avg. $L^2(\Omega_0)$-relative error & $1.5\cdot10^{-2}$ & $1.4\cdot10^{-3}$ & $1.7\cdot10^{-2}$\\
		\hline
	\end{tabular}
	\caption{Test case 2. Computational data related to POD-Galerkin-DEIM ROMs and Deep-HyROMnet, for $N=4$.}
	\label{tab:DB_hat_hyper-ROM_N4} 
\end{table}

\begin{figure}[t!]
	\centering
	\includegraphics[width=\textwidth]{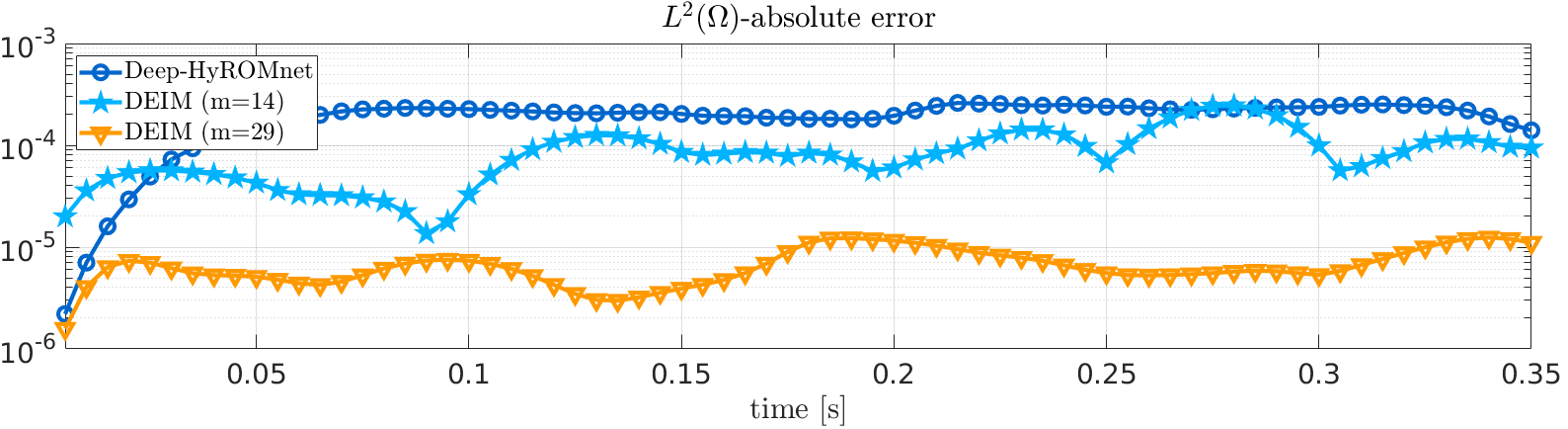}
	\caption{Test case 2. Evolution in time of the average $L^2(\Omega_0)$-absolute error for $N=4$ computed using POD-Galerkin-DEIM and Deep-HyROMnet.}
	\label{fig:DB_hat_hyper-ROM_N4_err_abs}
\end{figure}

\begin{figure}[t]
	\centering
	\includegraphics[width=\textwidth]{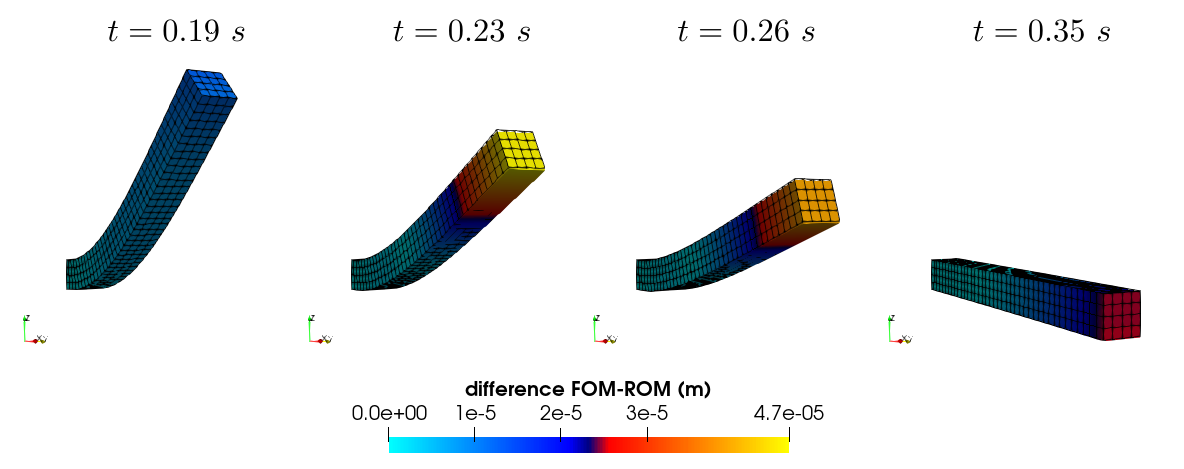}
	\caption{Test case 2. FOM (wireframe) and Deep-HyROMnet (colored) solutions computed at different times for $\bm\mu = [10.7375~\text{Pa}]$.}
	\label{fig:DB_hat_N4_DeepHyROMnet_mu36_error}
\end{figure}

The evolution in time of the average $L^2(\Omega_0)$-absolute error for DEIM and Deep-HyROMnet models is shown in Figure~\ref{fig:DB_hat_hyper-ROM_N4_err_abs}. 
The accuracy obtained using Deep-HyROMnet, although slightly lower than the ones achieved using a DEIM-based approximation, is satisfying in all the considered scenarios. Figure \ref{fig:DB_hat_N4_DeepHyROMnet_mu36_error} shows the FOM and the Deep-HyROMnet displacements at different time instances obtained for a given testing parameter.


\subsubsection{Test case 3: step function for the pressure load}

As last test case for the beam geometry, we consider a pressure load acting on the bottom surface area for only a third of the whole simulation time, that is
\begin{equation*}
\mathbf{g}(t;\bm\mu) = \widetilde{p}~\chi(t)_{\left(0,\frac{T}{3}\right]},
\end{equation*}
such that the resulting deformation features oscillations. This case is of particular interest in nonlinear elastodynamics, since the inertial term cannot be neglected, as it has a crucial impact on the deformation of the object. For the case at hand, we choose $t\in[0,0.27]$~s and a uniform time step $\Delta t = 3.6\cdot10^{-3}$~s, resulting in a total number of $75$ time iterations. For what concerns the input parameters, we vary the external load $\widetilde{p}\in[2,12]$ Pa and consider $G=10^4$ Pa and $K=5\cdot10^4$ Pa fixed.

\begin{figure}[b!]
\vspace{-0.2cm}
	\centering
	\begin{subfigure}{0.49\textwidth}
		\centering
		\includegraphics[width=0.95\textwidth]{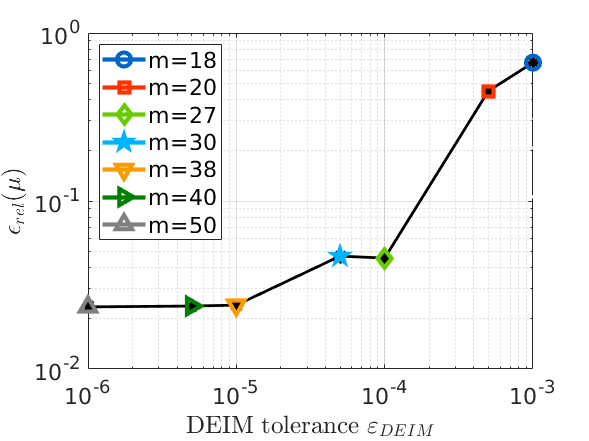}
	\end{subfigure}
	\begin{subfigure}{0.49\textwidth}
		\centering
		\includegraphics[width=0.95\textwidth]{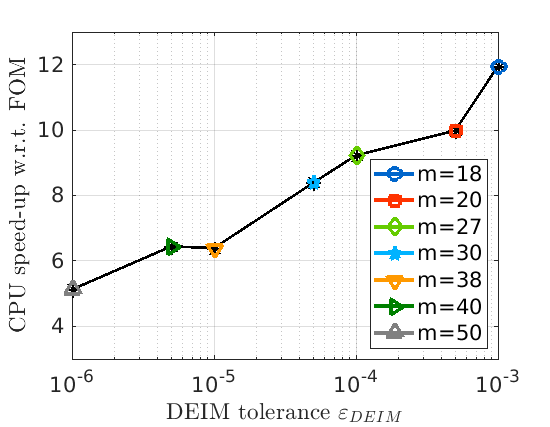}
	\end{subfigure}
	\caption{Test case 3. Average over 50 testing parameters of relative error $\epsilon_{rel}$ (left) and average speed-up (right) of POD-Galerkin-DEIM with $N=4$.}
	\label{fig:DB_step_DEIM_N4}
\end{figure}

We build the reduced basis $\mathbf{V}\in\mathbb{R}^{N_h\times N}$ from a training set of $n_s=50$ FOM solutions using $\varepsilon_{POD}=10^{-3}$, thus obtaining a reduced dimension of $N=4$, and perform POD-ROM simulations for a given set of $n_s'=300$ parameter samples to collect the nonlinear terms data necessary for the construction of both POD-Galerkin-DEIM and Deep-HyROMnet models. The DEIM basis $\mathbf{\bm\Phi}_{\mathcal{R}}$ for the approximation of the residual is computed by performing POD on the associated snapshots matrix with tolerance
\begin{equation*}
\varepsilon_{DEIM}\in\{10^{-3}, 5\cdot10^{-4}, 10^{-4}, 5\cdot10^{-5}, 10^{-5}, 5\cdot10^{-6}, 10^{-6}\},
\end{equation*}
where $\varepsilon_{DEIM}=10^{-3}$ is the larger POD tolerance that allows to guarantee the convergence of the reduced Newton algorithm for all testing parameters. The corresponding number basis for $\mathbf{R}$ is $m=18$, $20$, $27$, $30$, $38$, $40$, $50$, respectively. The results regarding the average relative error $\epsilon_{rel}$ and the computational speed-up, evaluated over 50 instances of the parameter, are shown in Figure \ref{fig:DB_step_DEIM_N4}. 

Like for the previous test cases, we compare POD-Galerkin-DEIM and Deep-HyROMnet ROMs, with respect to the displacement error and the CPU time ratio. 
\begin{table}[h]
	\centering
	\begin{tabular}{|l||c|c|c|} 
		\hline
		& DEIM ($m=18$) & DEIM ($m=38$) & Deep-HyROMnet \\
		\hline
		\hline
		Computational speed-up & $\times$12  & $\times$6 & $\times$1350\\
		\hline
		Avg. CPU time & 4~s & 8~s & 0.038~s\\
		\hline
		Time-avg. $L^2(\Omega_0)$-absolute error & $2.4\cdot10^{-3}$ & $1.3\cdot10^{-4}$ & $4.8\cdot10^{-4}$\\
		\hline
		Time-avg. $L^2(\Omega_0)$-relative error & $6.7\cdot10^{-1}$ & $2.4\cdot10^{-2}$ & $1.0\cdot10^{-1}$\\
		\hline
	\end{tabular}
	\caption{Test case 3. Computational data related to POD-Galerkin-DEIM ROMs and Deep-HyROMnet, for $N=4$.}
	\label{tab:DB_step_N4_hyper-ROM}
\end{table}
As reported in Table~\ref{tab:DB_step_N4_hyper-ROM}, Deep-HyROMnet outperforms DEIM substantially in terms of efficiency also in this case when handling the nonlinear terms. Indeed, Deep-HyROMnet yields a ROM that is more than $1000$ times faster than the FOM (this latter requiring $51$~s in average to be solved), still providing satisfactory results in terms of accuracy.

\begin{figure}[h!]
	\centering
	\begin{subfigure}{0.85\textwidth}
		\centering
		\includegraphics[width=\textwidth]{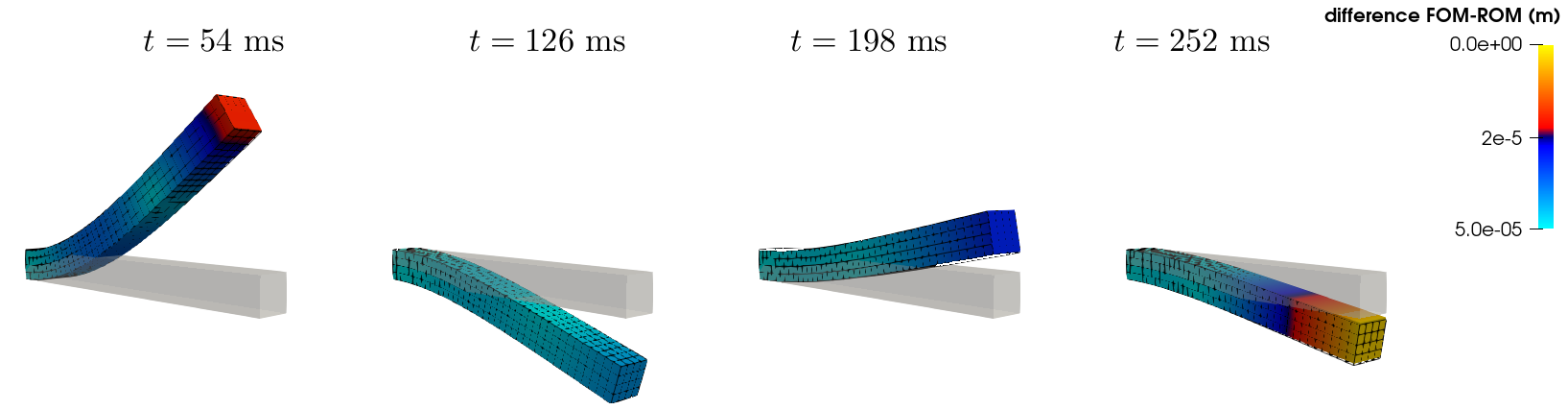}
		\caption{$\bm\mu = [5.1125~\text{Pa}]$}
	\end{subfigure}
	\begin{subfigure}{0.85\textwidth}
		\centering
		\includegraphics[width=\textwidth]{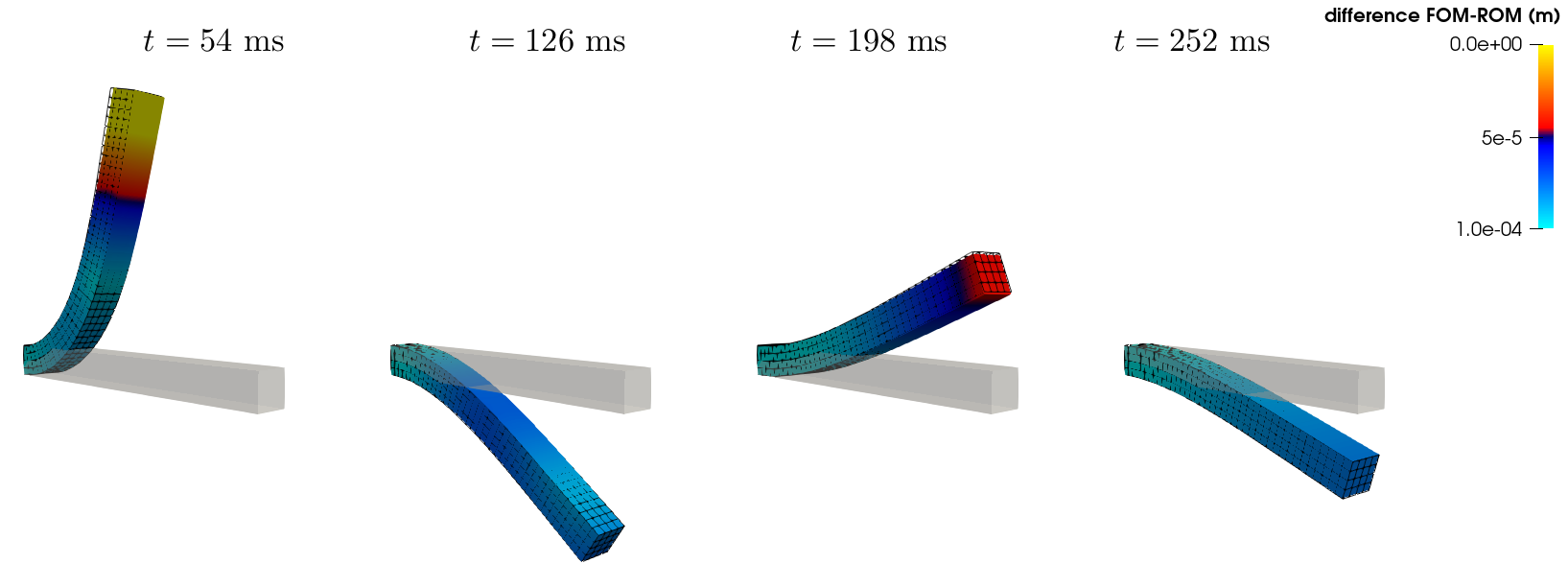}
		\caption{$\bm\mu = [11.4625~\text{Pa}]$}
	\end{subfigure}
	\caption{Test case 3. FOM (wireframe) and Deep-HyROMnet (colored) solutions computed at different times.}
	\label{fig:DB_step_N4_DeepHyROMnet_mu1_mu13}
\end{figure}

Figure~\ref{fig:DB_step_N4_DeepHyROMnet_mu1_mu13} represent the Deep-HyROMnet solution at different time instants for two different values of the parameter and show that the hyper-ROM is able to correctly capture the nonlinear behavior of the continuum body also when the inertial term cannot be neglected.


\subsection{Passive inflation and active contraction of an idealized left ventricle}

\begin{figure}
	\centering
	\includegraphics[width=0.725\textwidth]{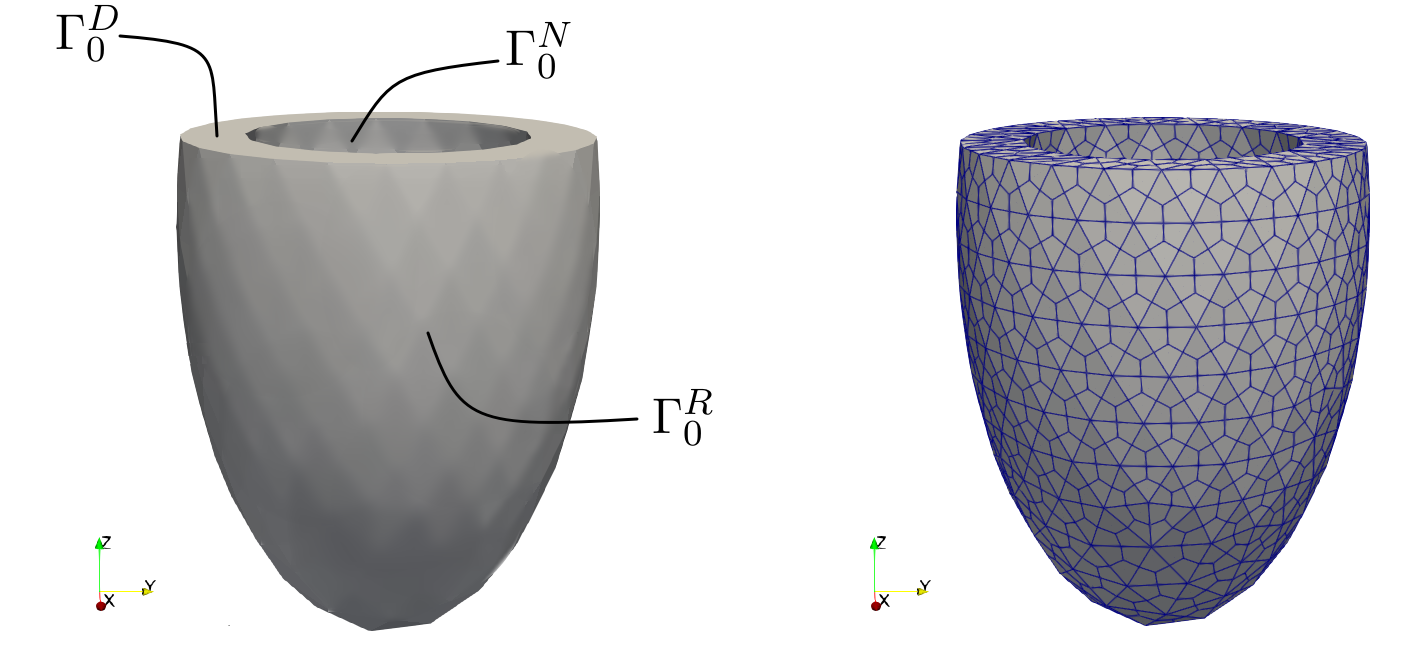}
	\caption{Passive inflation and active contraction of an idealized left ventricle. Idealized truncated ellipsoid geometry (left) and computational grid (right).} \label{fig:Prolate_BCs}
\end{figure}

The second problem we are interested in is the inflation and contraction of a prolate spheroid geometry representing an idealized left ventricle (see Figure~\ref{fig:Prolate_BCs}) where the boundaries $\Gamma_0^R$, $\Gamma_0^N$ and $\Gamma_0^D$ represent the epicardium, the endocardium and the base of a left ventricle, respectively, the latter being the artificial boundary resulting from truncation of the heart below the valves in a short axis plane.
We consider transversely isotropic material properties for the myocardial tissue, adopting a nearly-incompressible for\-mu\-la\-tion of the constitutive law proposed in \cite{guccione1995finite}, whose strain-energy density function is given by
\begin{equation*}
\mathcal{W}(\mathbf{F}) = \frac{C}{2}(e^{Q(\mathbf{F})} - 1),
\end{equation*}
with the following form for $Q$ to describe three-dimensional transverse isotropy with respect to the fiber coordinate system,
\begin{equation*}
Q = b_{f} E_{ff}^2 + b_{s} E_{ss}^2 + b_{n} E_{nn}^2 + b_{fs}(E_{fs}^2 + E_{sf}^2) + b_{fn}(E_{fn}^2 + E_{nf}^2) + b_{sn}(E_{sn}^2 + E_{ns}^2).
\end{equation*}
Here, $E_{ij}$, $i,j\in\{f,s,n\}$, are the components of the Green-Lagrange strain tensor (\ref{eq:strain}), the material constant $C>0$ scales the stresses and the coefficients $b_f$, $b_s$, $b_n$ are related to the material stiffness in the fiber, sheet and transverse directions, respectively. This leads to a (passive) first Piola-Kirchhoff stress tensor characterized by exponential nonlinearity. In order to enforce the incompressibility constraint, we consider an additional term $\mathcal{W}_{vol}(J)$ in the definition of the strain energy density function, which must grow as the deformation deviates from being isochoric. A common choice for $\mathcal{W}_{vol}$ is a convex function with null slope in $J=1$, e.g.,
\begin{equation*}
\mathcal{W}_{vol}(J) = \frac{K}{4}( (J-1)^2 + \ln^2(J) ),
\end{equation*} 
where the penalization factor is the bulk modulus $K>0$. Furthermore, to reproduce the typical twisting motion of the ventricular systole, we need to take into account a varying fiber distribution and contractile forces. The fiber direction is computed using the rule-based method proposed in \cite{rossi2014thermodynamically}, which depends on parameter angles $\bm\alpha^{epi}$ and $\bm\alpha^{endo}$. Active contraction is modeled through the active stress approach \cite{ambrosi2012active}, so that we add to the passive first Piola-Kirchoff stress tensor a time-dependent active tension, which is assumed to act only in the fiber direction
\begin{equation*}
\mathbf{P} = \left(\frac{\partial\mathcal{W}(\mathbf{F})}{\partial\mathbf{F}} + \frac{\partial\mathcal{W}_{vol}(J)}{\partial\mathbf{F}}\right) + T_a(t)(\mathbf{Ff}_0\otimes\mathbf{f}_0),
\end{equation*}
where $\mathbf{f}_0\in\mathbb{R}^3$ denotes the reference unit vector in the fiber direction and $T_a$ is a parametrized function that surrogates the active generation forces. In our case, since we are modeling only the systolic contraction, we define
\begin{equation*}
T_a(t) = \widetilde{T}_a~t/T, \quad t\in(0,T),
\end{equation*}
with $\widetilde{T}_a>0$. To model blood pressure inside the chamber we assume a linearly increasing external load 
\begin{equation*}
\mathbf{g}(t;\bm\mu) = \hat{p}~t/T, \quad t\in(0,T).
\end{equation*} 

Since we want to assess the performance of Deep-HyROMnet to reduce the myocardium contraction, we consider as unknown parameters those related to the active components of the strain energy function:
\begin{itemize}
	\item the maximum value of the active tension $\widetilde{T}_a\in[49.5\cdot10^3,70.5\cdot10^3]$ Pa, and
	\item the fiber angles $\bm\alpha^{epi}\in[-105.5,-74.5]^\circ$ and $\bm\alpha^{endo}\in[74.5,105.5]^\circ$.
\end{itemize}
All other parameters are fixed to the reference values taken from \cite{land2015verification}, namely $b_{f}=8$, $b_{s}=b_{n}=b_{sn}=2$, $b_{fs}=b_{fn}=4$, $C=2\cdot10^3$~Pa, $K=50\cdot10^3$~Pa and $\widetilde{p}=15\cdot10^3$. Regarding the time discretization, we choose $t\in[0,0.25]$~s and a uniform time step $\Delta t = 5\cdot10^{-3}$~s, resulting in a total number of $50$ time iterations. The FOM is built on a hexahedral mesh with $4804$ elements and $6455$ vertices, depicted in Figure \ref{fig:Prolate_BCs}, corresponding to a high-fidelity dimension $N_h=19365$, since $\mathbb{Q}_1$-FE (that is, linear FE on a hexahedral mesh) are used. In this case, the FOM requires almost $360$~s to compute the solution dynamics for each parameter instance.

\begin{figure}
	\centering
	\includegraphics[width=0.45\textwidth]{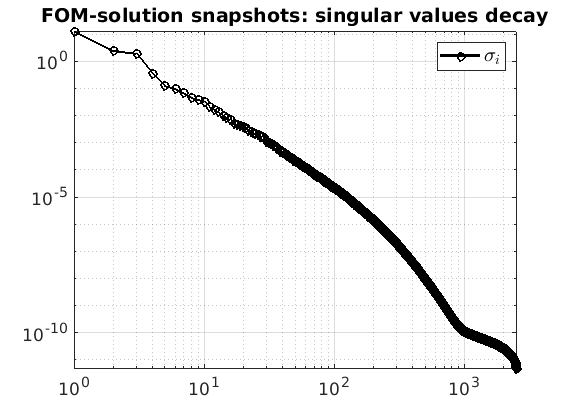}
	\caption{Passive inflation and active contraction of an idealized left ventricle. Decay of the singular values of the FOM solution snapshots matrix.}
	\label{fig:ProlateC_svd_uh} \vspace{-0.1cm}
\end{figure}

\begin{figure}[b!]
	\centering
	\begin{subfigure}{0.49\textwidth}
		\centering
		\includegraphics[width=0.95\textwidth]{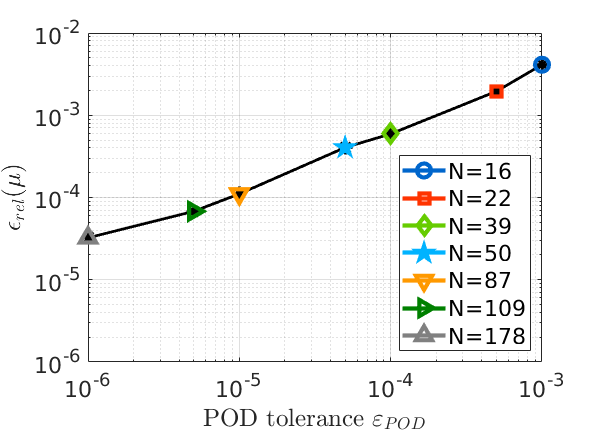}
	\end{subfigure}
	\begin{subfigure}{0.49\textwidth}
		\centering
		\includegraphics[width=0.95\textwidth]{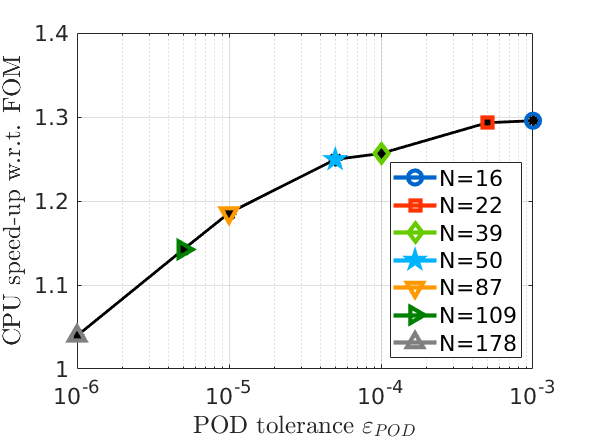}
	\end{subfigure}
	\caption{Passive inflation and active contraction of an idealized left ventricle. Average over 20 testing parameters of relative error $\epsilon_{rel}$ (left) and average speed-up (right) of ROM without hyper-reduction.}
	\label{fig:ProlateC_ROM}
\end{figure}

Given $n_s=50$ points obtained by sampling the parameter space $\mathcal{P}$, we construct the corresponding solution snapshots matrix $\mathbf{S}_u$ and compute the reduced basis $\mathbf{V}\in\mathbb{R}^{N_h\times N}$ using the POD method with tolerance
\begin{equation*}
\varepsilon_{POD}\in\{10^{-3}, 5\cdot10^{-4}, 10^{-4}, 5\cdot10^{-5}, 10^{-5}, 5\cdot10^{-6}, 10^{-6}\}.
\end{equation*}
From Figure~\ref{fig:ProlateC_svd_uh}, we observe a slower decay of the singular values of $\mathbf{S}_u$ with respect to the structural problems of Section~\ref{sec:DB}. In fact, we obtain larger reduced basis dimensions $N=16$, $22$, $39$, $50$, $87$, $109$ and $178$, respectively. 
The error and the CPU speed-ups averaged over a testing set of 20 parameters are both shown in Figure \ref{fig:ProlateC_ROM}, as functions of the POD tolerance $\varepsilon_{POD}$. 
As already discussed, the speed-up achieved by the ROM is negligible, since at each Newton iteration without hyper-reduction the ROM still depends on the high-fidelity dimension $N_h$. For what concerns the approximation error, we observe a reduction of almost two orders of magnitude when going from $N=16$ to $N=178$. \\

\begin{figure}[t!]
	\centering
	\includegraphics[width=0.45\textwidth]{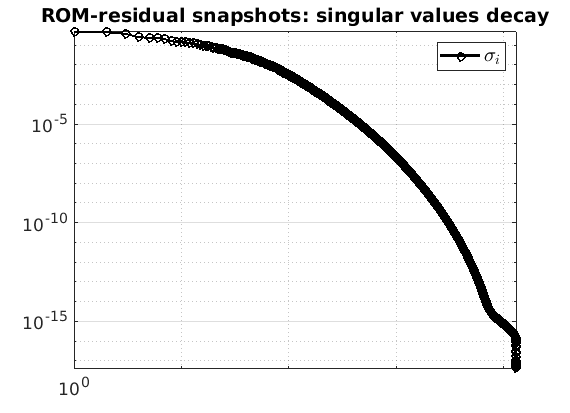}
	\caption{Passive inflation and active contraction of an idealized left ventricle. Decay of the singular values of the ROM residual snapshots matrix.}
	\label{fig:ProlateC_N16_svd_RN}
	\vspace{-0.25cm}
\end{figure}

\begin{figure}
	\centering
	\begin{subfigure}{0.49\textwidth}
		\centering
		\includegraphics[width=0.95\textwidth]{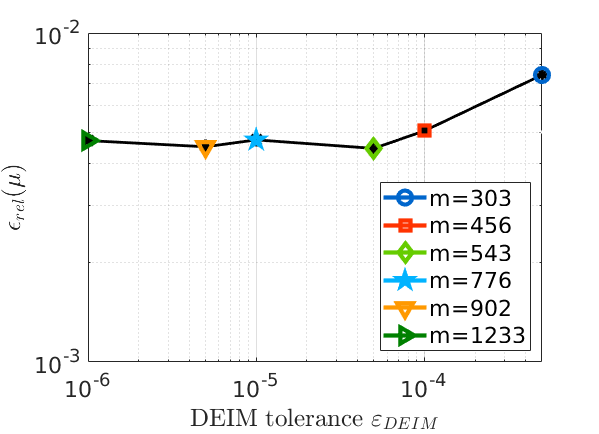}
	\end{subfigure}
	\begin{subfigure}{0.49\textwidth}
		\centering
		\includegraphics[width=0.95\textwidth]{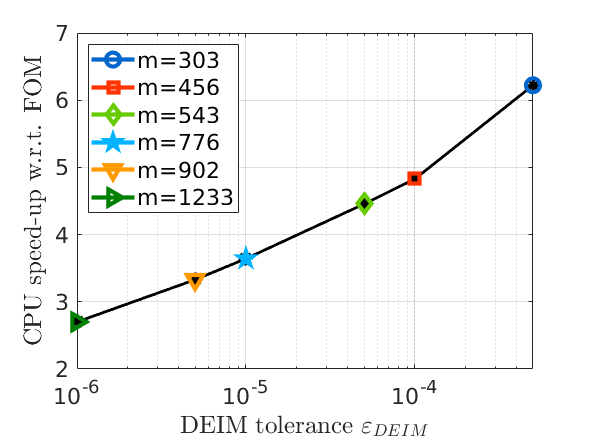}
	\end{subfigure}
	\caption{Passive inflation and active contraction of an idealized left ventricle. Average over 20 testing parameters of relative error $\epsilon_{rel}$ (left) and average speed-up (right) of POD-Galerkin-DEIM with $N=16$.}
	\label{fig:ProlateC_DEIM_N16}
	\vspace{-0.3cm}
\end{figure}

Given the reduced basis $\mathbf{V}\in\mathbb{R}^{N_h\times N}$ with $N=16$, we construct the POD-Galerkin-DEIM ap\-proxima\-tion by considering $n_s'=200$ parameter samples. Figure~\ref{fig:ProlateC_N16_svd_RN} shows the decay of the singular values of $\mathbf{S}_R$, that is, the snapshots matrix of the residual vectors $\mathbf{R}(\mathbf{Vu}_N^{n,(k)}(\bm\mu_{\ell'}),t^n;\bm\mu_{\ell'})$. We observe that the reported curve decreases very slowly, so that we expect that a large number of basis functions is required to correctly approximate the nonlinear operators.  
In fact, by computing $\bm\Phi_{\mathcal{R}}\in\mathbb{R}^{N_h\times m}$ using the following POD tolerances:
\begin{equation*}
\varepsilon_{DEIM}\in\{5\cdot10^{-4}, 10^{-4}, 5\cdot10^{-5}, 10^{-5}, 5\cdot10^{-6}, 10^{-6}\},
\end{equation*}
we obtain $m=303$, $456$, $543$, $776$, $902$ and $1233$, respectively. Higher values of $\varepsilon_{DEIM}$ (related to hopefully smaller dimensions $m$) were not sufficient to guarantee the convergence of the reduced Newton problem for all the parameter combinations considered. The average relative error over a set of 20 parameters and the computational speed-up are both reported in Figure~\ref{fig:ProlateC_DEIM_N16}. In particular, we observe that the relative error is between $4\cdot10^{-3}$ and $8\cdot10^{-3}$, as we could expect from the projection error reported in Figure~\ref{fig:ProlateC_ROM}, that is, POD-Galerkin-DEIM is able to achieve the same accuracy of the ROM without hyper-reduction.

The data reported in Table~\ref{tab:ProlateC} leads to the same conclusions regarding the computational bottleneck of the DEIM technique as those reported in Table~\ref{tab:DB_ramp}. In fact, assembling the residual on the reduced mesh requires around $85\%$ of the online CPU time, thus undermining the hyper-ROM efficiency. 
\begin{table}
	\centering
	\begin{tabular}{|l||c|c|c|} 
		\hline
		POD tolerance $\varepsilon_{DEIM}$ & $5\cdot10^{-4}$ & $5\cdot10^{-5}$ & $5\cdot10^{-6}$\\
		\hline
		DEIM interpolation dofs $m$ & $303$ & $543$ & $902$ \\
		\hline
		Reduced mesh elements (total: $4804$) & $914$ & $1345$ & $1855$ \\
		\hline
		\hline
		Online CPU time & $58$~s & $81$~s & $110$~s \\
		$\quad\circ$ system construction $[*]$ & $89\%$   & $93\%$   & $94\%$\\
		$\quad\circ$ system solution           & $0.01\%$ & $0.01\%$ & $0.01\%$\\
		\hline
		\hline
		$[*]$ System construction for each Newton iteration   & $0.4$~s  & $0.6$~s  & $0.9$~s \\
		$\quad\circ$ residual assembling                      & $94\%$   & $94\%$   & $94\%$\\
		$\quad\circ$ Jacobian computing through AD            & $0.24\%$ & $0.24\%$ & $0.26\%$\\
		\hline
		\hline
		Computational speed-up & $\times$6.2 & $\times$4.5 & $\times$3.3 \\
		\hline
		Time-averaged $L^2(\Omega_0)$-absolute error & $1\cdot10^{-3}$ & $7\cdot10^{-4}$ & $6\cdot10^{-4}$\\
		\hline
		Time-averaged $L^2(\Omega_0)$-relative error & $8\cdot10^{-3}$ & $5\cdot10^{-3}$ & $5\cdot10^{-3}$\\
		\hline
	\end{tabular}
\caption{Passive inflation and active contraction of an idealized left ventricle. Computational data related to POD-Galerkin-DEIM with $N=16$ and different values of $m$.}
\label{tab:ProlateC}
\end{table}
   
Finally, Table \ref{tab:ProlateC_N16_hyper-ROM} reports the computational data of POD-Galerkin-DEIM ROMs obtained for a number of magic points equals to $m=303$ and $m=543$, and of the Deep-HyROMnet, clearly showing that the latter outperforms the classical reduction strategy regarding the computational speed-up. 

\begin{figure}[b!]
	\centering
	\includegraphics[width=0.925\textwidth]{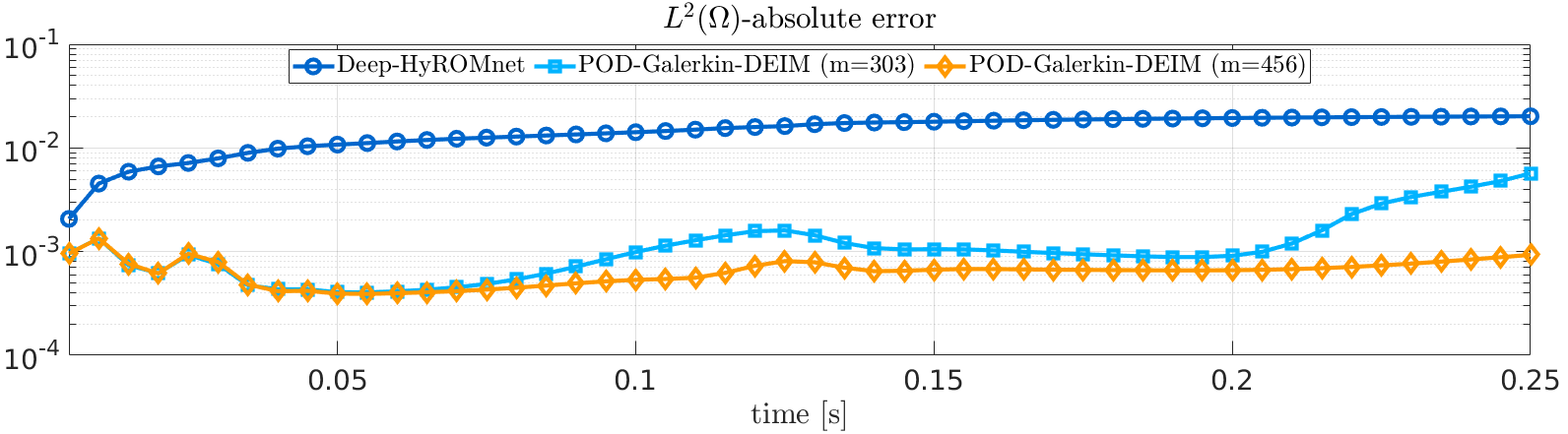}
	\caption{Passive inflation and active contraction of an idealized left ventricle. Evolution in time of the average $L^2(\Omega_0)$-absolute error computed using DEIM-based and DNN-based hyper-ROMs, for $N=16$.}
	\label{fig:ProlateC_N16_hyper-ROM_err_abs}
\end{figure}

\begin{figure}[t!]
	\centering
	\includegraphics[width=0.95\textwidth]{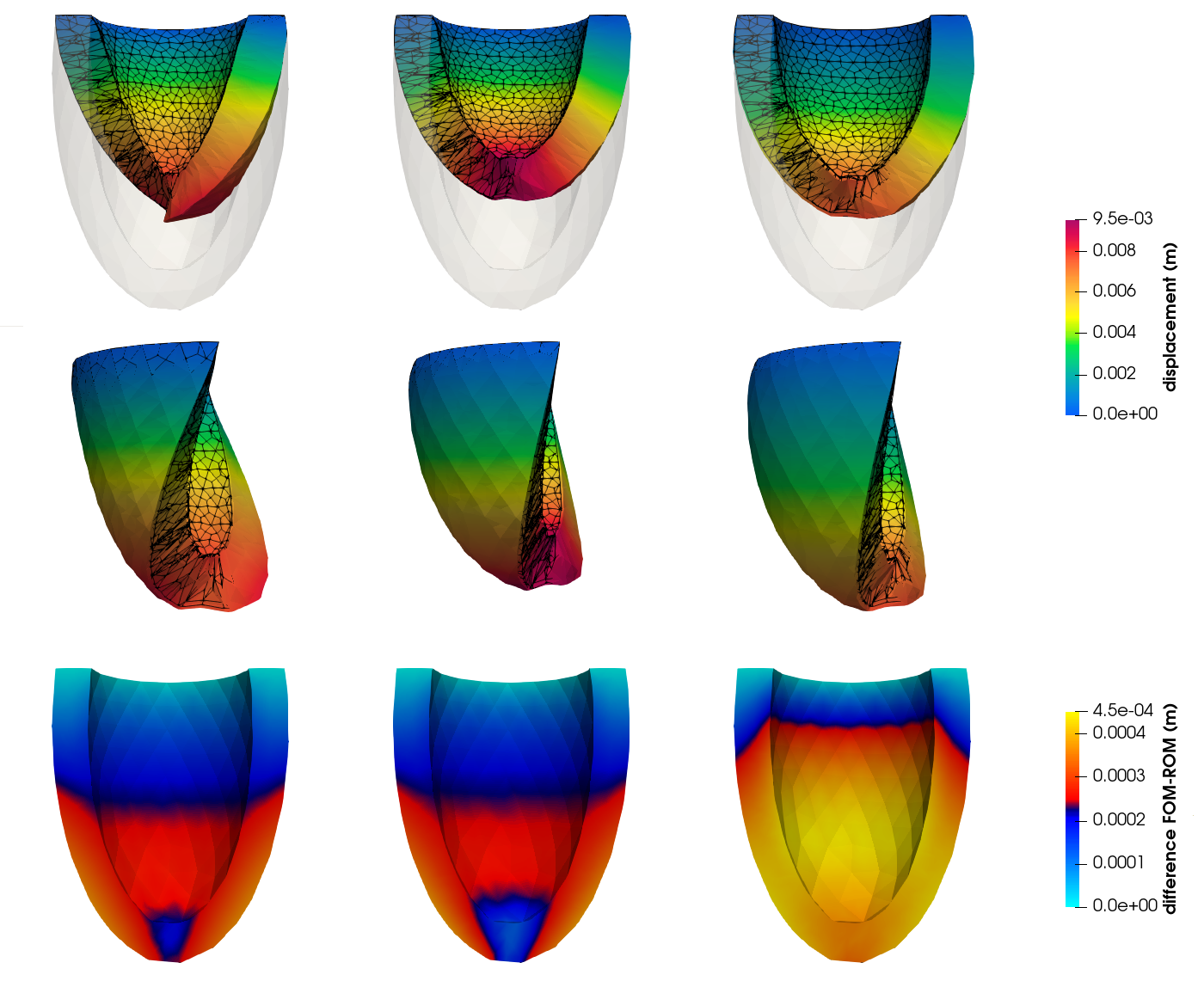}
	\caption{Passive inflation and active contraction of an idealized left ventricle. FOM (wireframe) and Deep-HyROMnet (colored) displacements (frontal view on top, lateral view in the middle) and corresponding difference (bottom) at time $T=0.25$~s for $\bm\mu = [61942.5~\text{Pa},-77.5225^\circ, 87.9075^\circ]$ (left), $\bm\mu = [59737.5~\text{Pa},-102.3225^\circ, 91.1625^\circ]$ (center) and $\bm\mu = [50497.5~\text{Pa},-100.9275^\circ, 80.0025^\circ]$ (right).}
	\label{fig:ProlateC_N16_DeepHyROMnet_mu1_mu10_mu19}
\end{figure}

In fact, Deep-HyROMnet is able to approximate the solution dynamics in $0.1$~s, that is even faster than real-time, while a POD-Galerkin-DEIM ROM requires $1$~min in average, where the final simulation time $T$ is set equal to $0.25$~s. Although the Deep-HyROMnet error is one order of magnitude higher than the one evaluated by a DEIM-based hyper-ROM (see Figure~\ref{fig:ProlateC_N16_hyper-ROM_err_abs}), the results are satisfactory in terms of accuracy. In Figures~\ref{fig:ProlateC_N16_DeepHyROMnet_mu1_mu10_mu19} the FOM and the DNN-based hyper-ROM displacements at time $T=0.25$~s are reported for three different values of the parameters, together with the error between the high-fidelity and the reduced solutions. 
\begin{table}[h!]
	\centering
	\begin{tabular}{|l||c|c|c|} 
		\hline
		& DEIM ($m=303$) & DEIM ($m=543$) & Deep-HyROMnet \\
		\hline
		\hline
		Speed-up  & $\times$6  & $\times$5 & $\times$3554 \\
		\hline
		Avg. CPU time & 58~s & 75~s & 0.1~s\\
		\hline
		mean$_{\bm\mu}$ $\epsilon_{abs}(\bm\mu)$ & $1.3\cdot10^{-3}$ & $6.6\cdot10^{-4}$ & $1.5\cdot10^{-2}$ \\
		\hline
		mean$_{\bm\mu}$ $\epsilon_{rel}(\bm\mu)$ & $7.5\cdot10^{-3}$ & $5.0\cdot10^{-3}$ & $6.4\cdot10^{-2}$ \\
		\hline
	\end{tabular}
	\caption{Passive inflation and active contraction of an idealized left ventricle. Computational data related to DEIM-based and DNN-based hyper-ROMs, for $N=16$.}
	\label{tab:ProlateC_N16_hyper-ROM} 
\end{table}


To conclude, we repeat that the approximation of the reduced nonlinear operators with Deep-HyROMnet does not depend directly on the high-fidelity dimension $N_h$, but rather on reduced basis dimension $N$. 
To test its performances using a higher FOM dimension, we address the solution to the problem described in this Section, however considering a finer hexahedral mesh with $9964$ elements and $13025$ vertices, thus obtaining $N_h=39075$ as FOM dimension. In this case, about 13 minutes are required to compute the high-fidelity dynamics. On the other hand, a reduced basis of dimension $N=16$ is computed for $\varepsilon_{POD}=10^{-3}$; the computational data, averaged over a testing set of $20$ parameter samples, are reported in Table \ref{tab:ProlateC_hyper-ROM_fine}. Almost unexpectedly, the online CPU time required by Deep-HyROMnet doubles as we double $N_h$. This may be due to the higher time required to perform matrix-vector multiplication for the reconstruction of the reduced solutions $\mathbf{Vu}_N^n(\bm\mu)$, for $n=1,\dots,N_t$. Further analysis should be performed to investigate this issue. Nonetheless, is it worth saying that the overall computational speed-up of Deep-HyROMnet increases as the FOM dimension $N_h$ grows, while the number $N$ of reduced basis function remains small, so that reduced solutions can be computed extremely fast. For what concerns the approximation accuracy of the hyper-ROMs with respect to the associate FOMs, we obtain almost the same results, showing that Deep-HyROMnet is able to deal with higher high-fidelity dimensions.
\begin{table}
	\centering
	\begin{tabular}{|l||c|c|} 
		\hline
		& \multicolumn{2}{c|}{Deep-HyROMnet} \\
		\hline
		\hline
		$N_h$ & $19365$ & $39075$ \\
		\hline
		FOM time & 5~min 54~s & 13~min 01~s  \\
		\hline
		$N$ & \multicolumn{2}{c|}{$16$} \\
		\hline
		Speed-up & $\times$3554 & $\times$3886 \\
		\hline
		Avg. CPU time & 0.1~s & 0.2~s \\
		\hline
		mean$_{\bm\mu}$ $\epsilon_{abs}(\bm\mu)$ & $1.5\cdot10^{-2}$ & $2.7\cdot10^{-2}$ \\
		\hline
		mean$_{\bm\mu}$ $\epsilon_{rel}(\bm\mu)$ & $6.4\cdot10^{-2}$ & $8.3\cdot10^{-2}$ \\
		\hline
	\end{tabular}
	\caption{Passive inflation and active contraction of an idealized left ventricle. Computational data related to Deep-HyROMnet for $N_h=19365$ and $N_h=39075$.}
	\label{tab:ProlateC_hyper-ROM_fine}
\end{table}


\section{Conclusions}\label{sec:conclusion}

In this work we have addressed the solution to the parametrized elastodynamics equation, correlated with nonlinear constitutive law, by means of a new projection-based reduced order model (ROM), developed to accurately capture the state solution dynamics at a reduced computational cost with respect to full-order models (FOMs) providing expensive high-fidelity approximations.

We focused on Galerkin-reduced basis (RB) methods, characterized by a projection of the differential problem onto a low-dimensional subspace built, e.g., by performing proper orthogonal decomposition (POD) on a set of FOM solutions, and by the splitting of the reduction procedure into a costly offline phase and an inexpensive online phase. Numerical experiments showed that, despite their highly nonlinear nature, elastodynamics problems can be reduced by exploiting projection-based strategies in an effective way, with POD-Galerkin ROMs achieving very good accuracy even in presence of a handful of basis functions. However, when dealing with nonlinear problems, a further level of approximation is required to make the online stage independent of the high-fidelity dimension. 

Hyper-reduction techniques, such as the discrete empirical interpolation method (DEIM), are necessary to efficiently handle the nonlinear operators. However, a serious issue is represented by the assembling (albeit onto a reduced mesh) of the approximated nonlinear operators in this framework. This observation suggested the idea of relying on surrogate models to perform operator approximation, overcoming the need to assemble the nonlinear terms onto the computational mesh.

Pursuing this strategy, we have proposed a new projection-based, deep learning-based ROM, \textit{Deep-HyROMnet}, which combines the Galerkin-RB approach with deep neural networks (DNNs) to assemble the reduced Newton system in an efficient way, thus avoiding the computational burden entailed by classical hyper-reduction strategies. This approach allows to rely on physics-based (thus, consistent) ROMs retaining the underlying structure of the physical model, as DNNs are employed only for the approximation of the reduced nonlinear operators, so that the problem displacement at each time instance is computed by solving the reduced nonlinear system. Regarding the offline cost of this hybrid reduction strategy, we point out that:
\begin{itemize}
	\item FOM solutions are required only for the construction of the reduced basis functions;
	\item since the nonlinear operators are collected during Newton iterations at each time step, a smaller number of ROM simulations with respect to purely data-driven approaches is sufficient for training the DNNs;
	\item being the training data low-dimensional, we can avoid the overwhelming training times and costs that would be required by the DNN if FOM arrays were used. 	
\end{itemize}

Deep-HyROMnet has been successfully applied in a nonlinear solid mechanics context, showing remarkable improvement in terms of online CPU time with respect to POD-Galerkin-DEIM ROMs. Our goal in future works is to apply the developed strategy to other classes of nonlinear problems for which traditional hyper-reduction techniques represent a computational bottleneck.

\bibliography{CFM_bibliography.bib}

\begin{thebibliography}{10}

\bibitem{benner2015survey}
P.~Benner, S.~Gugercin, and K.~Willcox, ``A survey of projection-based model
  reduction methods for parametric dynamical systems,'' {\em SIAM {R}eview},
  vol.~57, no.~4, pp.~483--531, 2015.

\bibitem{benner2017model}
P.~Benner, M.~Ohlberger, A.~Patera, G.~Rozza, and K.~Urban~(Eds.), {\em Model
  {R}eduction of {P}arametrized {S}ystems}.
\newblock Springer, 2017.

\bibitem{quarteroni2016reduced}
A.~Quarteroni, A.~Manzoni, and F.~Negri, {\em Reduced {B}asis {M} ethods for
  {P}artial {D}ifferential {E}quations. {A}n {I}ntroduction}.
\newblock 2016.

\bibitem{hesthaven2016certified}
J.~Hesthaven, G.~Rozza, and B.~Stamm, {\em Certified reduced basis methods for
  parametrized partial differential equations}.
\newblock Springer, 2016.

\bibitem{farhat20205}
C.~Farhat, S.~Grimberg, A.~Manzoni, and A.~Quarteroni, ``Computational
  bottlenecks for {PROMs}: Pre-computation and hyperreduction,'' in {\em {Model
  Order Reduction. Volume 2: Snapshot-Based Methods and Algorithms}}
  (P.~Benner, S.~Grivet-Talocia, A.~Quarteroni, G.~Rozza, W.~Schilders, and
  L.~Silveira, eds.), pp.~181--244, Berlin: De~Gruyter, 2020.

\bibitem{pinkus2012n}
A.~Pinkus, {\em n-Widths in Approximation Theory}.
\newblock Berlin-Heidelberg: Springer-Verlag, 1985.

\bibitem{fresca2020deep}
S.~Fresca, A.~Manzoni, L.~Dede', and A.~Quarteroni, ``Deep learning-based
  reduced order models in cardiac electrophysiology,'' {\em PloS one}, vol.~15,
  no.~10, p.~e0239416, 2020.

\bibitem{amsallem2012nonlinear}
D.~Amsallem, M.~Zahr, and C.~Farhat, ``Nonlinear model order reduction based on
  local reduced-order bases,'' {\em International Journal for Numerical Methods
  in Engineering}, vol.~92, no.~10, pp.~891--916, 2012.

\bibitem{pagani2018numerical}
S.~Pagani, A.~Manzoni, and A.~Quarteroni, ``Numerical approximation of
  parametrized problems in cardiac electrophysiology by a local reduced basis
  method,'' {\em Computer Methods in Applied Mechanics and Engineering},
  vol.~340, pp.~530--558, 2018.

\bibitem{vlachas2021local}
K.~Vlachas, K.~Tatsis, K.~Agathos, A.~Brink, and E.~Chatzi, ``A local basis
  approximation approach for nonlinear parametric model order reduction,'' {\em
  Journal of Sound and Vibration}, vol.~502, p.~116055, 2021.

\bibitem{lee2020model}
K.~Lee and K.~Carlberg, ``Model reduction of dynamical systems on nonlinear
  manifolds using deep convolutional autoencoders,'' {\em Journal of
  Computational Physics}, vol.~404, p.~108973, 2020.

\bibitem{kim2020efficient}
Y.~Kim, Y.~Choi, D.~Widemann, and T.~Zohdi, ``Efficient nonlinear manifold
  reduced order model,'' {\em arXiv preprint arXiv:2011.07727}, 2020.

\bibitem{fresca2021comprehensive}
S.~Fresca, L.~Dede', and A.~Manzoni, ``A comprehensive deep learning-based
  approach to reduced order modeling of nonlinear time-dependent parametrized
  {P}{D}{E}s,'' {\em Journal of Scientific Computing}, vol.~87, no.~2,
  pp.~1--36, 2021.

\bibitem{fresca2021pod}
S.~Fresca and A.~Manzoni, ``{P}{O}{D}-{D}{L}-{R}{O}{M}: enhancing deep
  learning-based reduced order models for nonlinear parametrized {P}{D}{E}s by
  proper orthogonal decomposition,'' {\em Computer Methods in Applied Mechanics
  and Engineering}, vol.~388, no.~114181, 2022.

\bibitem{franco2021deep}
N.~Franco, A.~Manzoni, and P.~Zunino, ``A deep learning approach to reduced
  order modelling of parameter dependent partial differential equations,'' {\em
  arXiv preprint arXiv:2103.06183}, 2021.

\bibitem{barrault2004empirical}
M.~Barrault, Y.~Maday, N.~Nguyen, and A.~Patera, ``An ‘empirical
  interpolation’ method: application to efficient reduced-basis
  discretization of partial differential equations,'' {\em Comptes Rendus
  Mathematique}, vol.~339, no.~9, pp.~667--672, 2004.

\bibitem{chaturantabut2010nonlinear}
S.~Chaturantabut and D.~Sorensen, ``Nonlinear model reduction via discrete
  empirical interpolation,'' {\em SIAM Journal on Scientific Computing},
  vol.~32, no.~5, pp.~2737--2764, 2010.

\bibitem{negri2015efficient}
F.~Negri, A.~Manzoni, and D.~Amsallem, ``Efficient model reduction of
  parametrized systems by matrix discrete empirical interpolation,'' {\em
  Journal of Computational Physics}, vol.~303, pp.~431--454, 2015.

\bibitem{astrid2008missing}
P.~Astrid, S.~Weiland, K.~Willcox, and T.~Backx, ``Missing point estimation in
  models described by proper orthogonal decomposition,'' {\em IEEE Transactions
  on Automatic Control}, vol.~53, no.~10, pp.~2237--2251, 2008.

\bibitem{carlberg2011efficient}
K.~Carlberg, C.~Bou-Mosleh, and C.~Farhat, ``Efficient non-linear model
  reduction via a least-squares {P}etrov--{G}alerkin projection and compressive
  tensor approximations,'' {\em International Journal of Numerical Methods in
  Engineering}, vol.~86, no.~2, pp.~155--181, 2011.

\bibitem{farhat2015structure}
C.~Farhat, T.~Chapman, and P.~Avery, ``Structure-preserving, stability, and
  accuracy properties of the energy-conserving sampling and weighting method
  for the hyper reduction of nonlinear finite element dynamic models,'' {\em
  Int. J. Numer. Meth. Engng.}, vol.~102, no.~5, pp.~1077--1110, 2015.

\bibitem{drohmann2012reduced}
M.~Drohmann, B.~Haasdonk, and M.~Ohlberger, ``Reduced basis approximation for
  nonlinear parametrized evolution equations based on empirical operator
  interpolation,'' {\em SIAM Journal on Scientific Computing}, vol.~34, no.~2,
  pp.~A937--A969, 2012.

\bibitem{tiso2013discrete}
P.~Tiso and D.~Rixen, ``Discrete empirical interpolation method for finite
  element structural dynamics,'' in {\em Topics in Nonlinear Dynamics, Volume
  1}, pp.~203--212, Springer, 2013.

\bibitem{radermacher2016pod}
A.~Radermacher and S.~Reese, ``{P}{O}{D}-based model reduction with empirical
  interpolation applied to nonlinear elasticity,'' {\em International Journal
  for Numerical Methods in Engineering}, vol.~107, no.~6, pp.~477--495, 2016.

\bibitem{ghavamian2017pod}
F.~Ghavamian, P.~Tiso, and A.~Simone, ``{P}{O}{D}--{D}{E}{I}{M} model order
  reduction for strain-softening viscoplasticity,'' {\em Computer Methods in
  Applied Mechanics and Engineering}, vol.~317, pp.~458--479, 2017.

\bibitem{bonomi2017reduced}
D.~Bonomi, A.~Manzoni, and A.~Quarteroni, ``A matrix deim technique for model
  reduction of nonlinear parametrized problems in cardiac mechanics,'' {\em
  Computer Methods in Applied Mechanics and Engineering}, vol.~324,
  pp.~300--326, 2017.

\bibitem{cicci2021cardiacDEIM}
L.~Cicci, S.~Fresca, S.~Pagani, A.~Manzoni, and A.~Quarteroni,
  ``Projection-based reduced order models for parameterized nonlinear
  time-dependent problems arising in cardiac mechanics,'' tech. rep., 2021.
\newblock submitted. MOX Report N. 75/2021.

\bibitem{tiso2013modified}
P.~Tiso, R.~Dedden, and D.~Rixen, ``A modified discrete empirical interpolation
  method for reducing non-linear structural finite element models,'' in {\em
  International Design Engineering Technical Conferences and Computers and
  Information in Engineering Conference}, vol.~55973, p.~V07BT10A043, American
  Society of Mechanical Engineers, 2013.

\bibitem{peherstorfer2014localized}
B.~Peherstorfer, D.~Butnaru, K.~Willcox, and H.~Bungartz, ``Localized discrete
  empirical interpolation method,'' {\em SIAM Journal on Scientific Computing},
  vol.~36, no.~1, pp.~A168--A192, 2014.

\bibitem{hesthaven2018non}
J.~Hesthaven and S.~Ubbiali, ``Non-intrusive reduced order modeling of
  nonlinear problems using neural networks,'' {\em Journal of Computational
  Physics}, vol.~363, pp.~55--78, 2018.

\bibitem{guo2018reduced}
M.~Guo and J.~Hesthaven, ``Reduced order modeling for nonlinear structural
  analysis using gaussian process regression,'' {\em Computer Methods in
  Applied Mechanics and Engineering}, vol.~341, pp.~807--826, 2018.

\bibitem{guo2019data}
M.~Guo and J.~S. Hesthaven, ``Data-driven reduced order modeling for
  time-dependent problems,'' {\em Computer Methods in Applied Mechanics and
  Engineering}, vol.~345, pp.~75--99, 2019.

\bibitem{swischuk2019projection}
R.~Swischuk, L.~Mainini, B.~Peherstorfer, and K.~Willcox, ``Projection-based
  model reduction: Formulations for physics-based machine learning,'' {\em
  Computers \& Fluids}, vol.~179, pp.~704--717, 2019.

\bibitem{gao2020non}
H.~Gao, J.~Wang, and M.~Zahr, ``Non-intrusive model reduction of large-scale,
  nonlinear dynamical systems using deep learning,'' {\em Physica D: Nonlinear
  Phenomena}, vol.~412, p.~132614, 2020.

\bibitem{lu2021learning}
L.~Lu, P.~Jin, G.~Pang, Z.~Zhang, and G.~Karniadakis, ``Learning nonlinear
  operators via deeponet based on the universal approximation theorem of
  operators,'' {\em Nature Machine Intelligence}, vol.~3, no.~3, pp.~218--229,
  2021.

\bibitem{chen1995universal}
T.~Chen and H.~Chen, ``Universal approximation to nonlinear operators by neural
  networks with arbitrary activation functions and its application to dynamical
  systems,'' {\em IEEE Transactions on Neural Networks}, vol.~6, no.~4,
  pp.~911--917, 1995.

\bibitem{wang2021learning}
S.~Wang, H.~Wang, and P.~Perdikaris, ``Learning the solution operator of
  parametric partial differential equations with physics-informed deeponets,''
  {\em Science Advances}, vol.~7, no.~40, p.~eabi8605, 2021.

\bibitem{peherstorfer2016data}
B.~Peherstorfer and K.~Willcox, ``Data-driven operator inference for
  nonintrusive projection-based model reduction,'' {\em Computer Methods in
  Applied Mechanics and Engineering}, vol.~306, pp.~196--215, 2016.

\bibitem{benner2020operator}
P.~Benner, P.~Goyal, B.~Kramer, B.~Peherstorfer, and K.~Willcox, ``Operator
  inference for non-intrusive model reduction of systems with non-polynomial
  nonlinear terms,'' {\em Computer Methods in Applied Mechanics and
  Engineering}, vol.~372, p.~113433, 2020.

\bibitem{qian2019transform}
E.~Qian, B.~Kramer, A.~Marques, and K.~Willcox, ``Transform \& learn: A
  data-driven approach to nonlinear model reduction,'' in {\em AIAA Aviation
  2019 Forum}, p.~3707, 2019.

\bibitem{bai2021non}
Z.~Bai and L.~Peng, ``Non-intrusive nonlinear model reduction via machine
  learning approximations to low-dimensional operators,'' {\em Advanced
  Modeling and Simulation in Engineering Sciences}, vol.~8, no.~28, 2021.

\bibitem{bhattacharya2021model}
K.~Bhattacharya, B.~Hosseini, N.~Kovachki, and A.~Stuart, ``Model reduction and
  neural networks for parametric {P}{D}{E}s,'' {\em The SMAI Journal of
  Computational Mathematics}, vol.~7, pp.~121--157, 2021.

\bibitem{quarteroni2013numerical}
A.~Quarteroni, {\em Numerical {M}odels for {D}ifferential {P}roblems}.
\newblock Springer, 2nd~ed., 2013.

\bibitem{POD}
A.~Chatterjee, ``An introduction to the proper orthogonal decomposition,'' {\em
  Current Science}, vol.~78, no.~7, pp.~808--817, 2000.

\bibitem{halko2011finding}
N.~Halko, P.~Martinsson, and J.~Tropp, ``Finding structure with randomness:
  Probabilistic algorithms for constructing approximate matrix
  decompositions,'' {\em SIAM {R}eview}, vol.~53, no.~2, pp.~217--288, 2011.

\bibitem{grepl2007efficient}
M.~Grepl, Y.~Maday, N.~Nguyen, and A.~Patera, ``Efficient reduced-basis
  treatment of nonaffine and nonlinear partial differential equations,'' {\em
  ESAIM: Mathematical Modelling and Numerical Analysis}, vol.~41, no.~3,
  pp.~575--605, 2007.

\bibitem{gobat2021reduced}
G.~Gobat, A.~Opreni, S.~Fresca, A.~Manzoni, and A.~Frangi, ``Reduced order
  modeling of nonlinear microstructures through proper orthogonal
  decomposition,'' {\em Mechanical Systems and Signal Processing}, p.~accepted
  for publication, 2022.

\bibitem{manzoni2018reduced}
A.~Manzoni, D.~Bonomi, and A.~Quarteroni, ``Reduced order modeling for cardiac
  electrophysiology and mechanics: New methodologies, challenges and
  perspectives,'' in {\em Mathematical and Numerical Modeling of the
  Cardiovascular System and Applications} (D.~Boffi, L.~Pavarino, G.~Rozza,
  S.~Scacchi, and C.~Vergara, eds.), vol.~16 of {\em SEMA SIMAI Springer
  Series}, pp.~115--166, Springer, Cham, 2018.

\bibitem{Broyden}
C.~Broyden, ``A class of methods for solving nonlinear simultaneous
  equations,'' {\em Mathematics of Computation}, vol.~19, no.~92, pp.~577--593,
  1965.

\bibitem{Goodfellow-et-al-2016}
I.~Goodfellow, Y.~Bengio, and A.~Courville, {\em Deep Learning}.
\newblock MIT Press, 2016.

\bibitem{dealII92}
D.~Arndt, W.~Bangerth, B.~Blais, T.~Clevenger, M.~Fehling, A.~Grayver,
  T.~Heister, L.~Heltai, M.~Kronbichler, M.~Maier, P.~Munch, J.~Pelteret,
  R.~Rastak, I.~Thomas, B.~Turcksin, Z.~Wang, and D.~Wells, ``The
  \texttt{deal.II} library, version 9.2,'' {\em Journal of Numerical
  Mathematics}, vol.~28, no.~3, pp.~131--146, 2020.

\bibitem{guccione1995finite}
J.~Guccione, K.~Costa, and A.~McCulloch, ``Finite element stress analysis of
  left ventricular mechanics in the beating dog heart,'' {\em Journal of
  biomechanics}, vol.~28, no.~10, pp.~1167--1177, 1995.

\bibitem{rossi2014thermodynamically}
S.~Rossi, T.~Lassila, R.~Ruiz-Baier, A.~Sequeira, and A.~Quarteroni,
  ``Thermodynamically consistent orthotropic activation model capturing
  ventricular systolic wall thickening in cardiac electromechanics,'' {\em
  European Journal of Mechanics-A/Solids}, vol.~48, pp.~129--142, 2014.

\bibitem{ambrosi2012active}
D.~Ambrosi and S.~Pezzuto, ``Active stress vs. active strain in mechanobiology:
  constitutive issues,'' {\em Journal of Elasticity}, vol.~107, no.~2,
  pp.~199--212, 2012.

\bibitem{land2015verification}
S.~Land, V.~Gurev, S.~Arens, C.~Augustin, L.~Baron, R.~Blake, C.~Bradley,
  S.~Castro, A.~Crozier, M.~Favino, T.~Fastl, T.~Fritz, H.~Gao, A.~Gizzi,
  B.~Griffith, D.~Hurtado, R.~Krause, X.~Luo, M.~Nash, S.~Pezzuto, G.~Plank,
  S.~Rossi, D.~Ruprecht, G.~Seemann, N.~Smith, J.~Sundnes, J.~Rice,
  N.~Trayanova, D.~Wang, Z.~Wang, and S.~Niederer, ``Verification of cardiac
  mechanics software: benchmark problems and solutions for testing active and
  passive material behaviour,'' {\em Proceedings of the Royal Society A:
  Mathematical, Physical and Engineering Sciences}, vol.~471, no.~2184,
  p.~20150641, 2015.

\end{thebibliography}
\bibliographystyle{ieeetr}

\end{document}